\documentclass[11pt]
{article}
\usepackage{amsfonts}
\usepackage{amsmath}
\usepackage{amssymb}

\RequirePackage{srcltx}

\newtheorem{thr}{Theorem}
\newtheorem{lemma}{Lemma}[section]
\newtheorem{con}{Corollary}

\allowdisplaybreaks
\begin{document}

\begin{center}\Large{
Embedding theorems in constructive approximation}
\end{center}

\begin{center}
\bf{B. V. Simonov, S. Yu. Tikhonov
}
\end{center}

\vspace{10mm}
\begin{abstract}
In this paper necessary and sufficient conditions for the accuracy of embedding theorems of different
function classes are obtained.
The main result of the paper is the criterion for embedding between the generalized Weyl-Nikol'skii and Lipschitz classes.
To define the Weyl-Nikol'skii classes, we use the concept of a $(\lambda,\beta)$-derivative, which is a
generalization of the derivative in the sense of Weyl. As corollaries, we obtain estimates of norms
and moduli of smoothness of transformed Fourier series.
\end{abstract}

\vspace{10mm}



\section{Introduction}
 {\sf History of the question.} One of the main problems  of  constructive
theory of functions\footnote{The concept became well known through S.N. Bersntein's paper \cite[V. 2,
p.295-300; p. 349-360]{bern1}.} is in finding a relationship between differential properties of functions
 and its structural or construc\-ti\-ve characteristics.
 This topic started to develop
more than a century ago and in many cases the research was conducted as follows: authors
considered a given functional class and by investigating the properties of its elements obtained
embedding theorems with other functional classes.
We recommend the recent articles \cite{pin}, \cite{tel-ob}, \cite{juk} for a historical survey.

Below we write the three results, which influenced substantially further research and actually
gave rise  to the development of new areas within the
approximation theory.
\begin{align*}
&{\rm (A)}&
 \qquad f^{(r)}\in \textnormal{Lip} \,\alpha   \quad
&\Longleftrightarrow
 \qquad
E_n(f)= O\left(\frac{1}{n^{r+\alpha}}\right)& \qquad
&(0<\alpha<1, r\in \mathbf{Z}_+),\\
&{\rm (B)}&
 \qquad f^{(r)}\in \textnormal{Lip} \,\alpha
 \quad
&\Longleftrightarrow
 \qquad
\omega_{r+1}\Bigl(f,\frac{1}{n}\Bigr) = O
\left(\frac{1}{n^{r+\alpha}}\right)& \qquad &(0<\alpha<1, r\in
\mathbf{Z}_+),
\\
&{\rm (C)}& \qquad f\in \textnormal{Lip}\, \alpha
 \quad
&\,\,\Longrightarrow
 \qquad
\widetilde{f}\in \textnormal{Lip} \,\alpha& \qquad &(0<\alpha<1).
\end{align*}

{\it Criterion}
(A) was proved by
\footnote{
See \cite{bern3} for a detailed review of the question before
the 30ths of the 20th century.} D. Jackson (1911, \cite{djek}) and S.B. Stechkin (\cite{stechkin-jac})
in the necessity part, and by S. N. Bernstein (1912, \cite{bern},
\cite{bern2}) and Ch. de la Vall{\'e}e-Poussin (1919, \cite{vallee}) in the sufficiency part.

The theorems of this type are called direct and inverse theorems of
appro\-xi\-ma\-ti\-on theory.
Direct theorems for $L_p, 1\le p \le \infty$ (see e.g. \cite{devore}, \cite{tel-ob}, \cite{juk})
are written as follows:
\begin{eqnarray}
E_n(f)_p&\le&  C(k)\, \,\omega_{k} \Bigl(f,\frac{1}{n}\Bigr)_p,
 \qquad\,\,\,\, k, n\in \mathbf{N},
\label{a1}\\
E_n(f)_p&\le&   \frac{C(k)}{n^r} \,\, \omega_{k}
\Bigl(f^{(r)},\frac{1}{n} \Bigr)_p,  \quad\,\, k, n, r \in
\mathbf{N}. \label{a2}
\end{eqnarray}
Inverse theorems for $L_p, 1\le p \le \infty$ (see e.g. \cite{devore}, \cite{tel-ob}, \cite{juk}):
\begin{equation} \label{a3}
\omega_{k} \Bigl(f,\frac{1}{n}\Bigr)_p \le  \frac{C(k)}{n^k}
\sum\limits_{\nu=0}^n (\nu+1)^{k-1} E_\nu(f)_p ,
 \quad k, n\in \mathbf{N},
\end{equation}
\begin{equation} \label{a4}
\omega_{k} \Bigl(f^{(r)},\frac{1}{n}\Bigr)_p \le C(k) \left(
\frac{1}{n^k} \sum\limits_{\nu=0}^n (\nu+1)^{k+r-1} E_\nu(f)_p +
\sum\limits_{\nu=n+1}^\infty \nu^{r-1} E_\nu(f)_p \right),  \,\, k,
n\in \mathbf{N}.
\end{equation}
Here and further, the best trigonometric  approximation $E_n(f)_p$ and the modulus of smoothness
$\omega_k(f,\delta)_p$ are  defined as follows:
$$
E_n(f)_p=\min\left( \left\|f-T \right\|_p : T\in \mathbf{T}_n
\right), \qquad \mathbf{T}_n=span\left\{\cos  mx, \sin mx :
|m|\le n  \right\}
$$
and
\begin{equation} \label{mod}
\omega_k\left(f,\delta\right)_p = \sup \limits_{|h|\le\delta}
\left\| \triangle^k_h f (x) \right\|_p,
\end{equation}
$$ \qquad \triangle^k_h f (x) = \triangle^{k-1}_h \left(
\triangle_h f (x)\right) \qquad\mbox{ш} \qquad \triangle_h f
(x)=f(x+h)-f(x),
$$
respectively.

In the case of $1<p<\infty$ one can write (\cite[p. 210]{devore}, \cite{tel-ob}, \cite{timan}) the following improvement of estimates (\ref{a1}) and (\ref{a3}):
\begin{eqnarray} \label{asharp}
\frac{C(k)}{n^k} \left( \sum\limits_{\nu=0}^n (\nu+1)^{k\tau-1}
E_\nu(f)_p^\tau \right)^\frac{1}{\tau}&\le&
 \omega_{k}
\Bigl(f,\frac{1}{n}\Bigr)_p\\&\le&  \frac{C(k)}{n^k} \left(
\sum\limits_{\nu=0}^n (\nu+1)^{k\theta-1} E_\nu(f)_p^\theta
\right)^\frac{1}{\theta}
\nonumber
\end{eqnarray}
where $k, n\in \mathbf{N}$, $\theta=\min(2,p)$, and $\tau=\max(2,p)$.

In view of estimates (\ref{a2}) and (\ref{a4}), we note that investigation of the question
on existence of the $r$-th derivative of $f$ from a given function space has been
initiated by Bernstein \cite{bern}. He proved that the condition $\sum\limits_{\nu=1}^{\infty}
\nu^{r-1} E_\nu(f)_\infty<\infty$ implies $f^{(r)}\in C$.
Later on, for $L_p\,(1\le p \le \infty)$,
the following results were obtained (see the review \cite{juk} and the paper by O.V. Besov
\cite{bes1}). For convenience, we write these embeddings  in terms of the Besov space $B_{p,\theta}^r$
and the Sobolev space $W_p^r$:
\begin{align*}
\qquad &B_{p,1}^r\subset W^r_p \subset B_{p,\infty}^r\qquad
p=1,\infty, \qquad
 \\
\qquad &B_{p,p}^r\subset W^r_p \subset B_{p,2}^r \qquad 1<p\le 2,
\qquad
\\
\qquad &B_{p,2}^r\subset W^r_p \subset B_{p,p}^r \qquad 2\le
p<\infty.
\end{align*}

{\it Criterion}
(B) was proved by
A. Zygmund (1945, \cite{zygmu-smooth}).
He was one of the first to use the modulus of smoothness concept of an integer order introduced by  Bernstein in 1912
(\cite{bern}).
At present, the moduli of smoothness properties are well-studied (\cite{joh}, \cite{juk})
 and the result (B) follows from the following inequalities (see \cite[Chapters 2 and 6]{devore},
\cite{johnen}):
($1\le p\le\infty$)
\begin{equation} \label{b1}
\omega_{k+r} \Bigl(f,\frac{1}{n}\Bigr)_p \le  \frac{C(k,r)}{n^r}
\omega_{k} \Bigl(f^{(r)},\frac{1}{n}\Bigr)_p,
 \qquad k, r, n \in \mathbf{N}
\end{equation}
\begin{equation} \label{b2}
\omega_{k} \Bigl(f^{(r)},\frac{1}{n}\Bigr)_p \le  C(k,r)
\sum\limits_{\nu=n+1}^\infty \nu^{r-1} \omega_{k+r}
\Bigl(f,\frac{1}{\nu}\Bigr)_p,
 \,\, k, r, n \in \mathbf{N}.
\end{equation}
Comparing the last two inequalities and inequalities (\ref{a2}) and (\ref{a4}) we see that from
(\ref{b1}) and (\ref{b2}), using (\ref{a1}) and (\ref{a3}), it is easy to get (\ref{a2}) and
(\ref{a4}).

We also mention  the paper by J. Marcinkiewicz (1938, \cite{ma}), where the following two inequalities were proved:
\begin{equation} \label{marc1}
 \|f'\|_p\le C(p) \left(\int\limits_0^{2\pi}\frac{ \omega_{2}
(f,u)_p^p}{u^{p+1}}{du}\right)^\frac 1p, \qquad 1<p\le 2
\end{equation}
and
\begin{equation} \label{marc2}
 \left(\int\limits_0^{2\pi}\frac{ \omega_{2}
(f,u)_p^p}{u^{p+1}}{du}\right)^\frac 1p \le C(p)
 \|f'\|_p, \qquad 2\le p<\infty.
\end{equation}
It is easy to show \cite[Rem 3.5]{studia} that (\ref{b2}) and
(\ref{marc1}) are corollaries of the estimate  \cite[Th 3.1]{studia}
\begin{equation} \label{b22}
\omega_{k} \Bigl(f^{(r)},\frac{1}{n}\Bigr)_p \le  C(k,r)
\left(\sum\limits_{\nu=n+1}^\infty \nu^{r\theta-1}
\omega_{k+r}^\theta
\Bigl(f,\frac{1}{\nu}\Bigr)_p\right)^\frac{1}{\theta},
 \,\,\,\, \end{equation}
 where $ k, r, n \in \mathbf{N}$,
$\theta=\min(2,p)$, and $1\le p<\infty$.

{\it Criterion} (C) was proved by I.I. Privalov (1919, \cite{priv}).
The following inequality, which implies embedding (C), was obtained by A. Zygmund (\cite{zygmu-smooth})
and N.K. Bary and S.B. Stechkin (\cite{ba-st}) ($p=1,\infty$):
\begin{equation} \label{c1}
\omega_{k} \Bigl(\widetilde{f}^{(r)},\frac{1}{n}\Bigr)_p \le
C(k, r) \left( n^{-k} \sum\limits_{\nu=1}^n
\nu^{k+r-1}\omega_{k}\Bigl(f,\frac{1}{\nu}\Bigr)_p +
\sum\limits_{\nu=n+1}^\infty \nu^{r-1} \omega_{k}
\Bigl(f,\frac{1}{\nu}\Bigr)_p \right),
 \,\, k, n \in \mathbf{N}, r \in \mathbf{Z_+},
\end{equation}
Here and further, $\widetilde{f}$ denotes the conjugate function to  $f$
\cite[V. 1, Ch. 2]{zygm}.

Next, we note the paper by G.H. Hardy and J.S. Littlewood (1928, \cite{h-l}) in which seemingly, for
the first time, some problems of constructive approximation theory were formulated and solved in terms of
embedding theorems.
 Some historical aspects of this approach were presented in the paper \cite{og}.

Finally, we mention some improvements of above written inequalities
for the case of the generalized derivatives and moduli of smoothness, as these estimates  are of a particular interest
 for us: (\ref{a1}) and
(\ref{a3}) are proved for any  $k>0$ in \cite{butz} and \cite{tab};
analogues of inequality (\ref{a2}) for the $f^\psi_\beta$-derivatives
are  shown in \cite[6.3]{step}, and analogues of (\ref{a4}) and (\ref{c1}) are proved
in \cite{samko}, \cite{step-1995}-\cite{zhukina}, \cite{taberski2}.
\\
{\sf Embedding theorems for functional classes.}
The results (A) - (C) as well as their
generalizations mentioned above can be written as the embedding theorems of the following
functional classes:
\begin{align*}
W^r_p &=\Bigl\{ f\in L_p: f^{(r)}\in L_p \Bigr\}, \\
\widetilde{W}^r_p &=\left\{ f\in L_p:\widetilde{f}^{(r)}\in L_p \right\}, \\
{W}^r_pH_\alpha[\varphi] &=\Bigl\{ f\in W^r_p:
\omega_\alpha \left(f^{(r)},\delta \right)_p = O[\varphi(\delta)] \Bigr\}, \\
\widetilde{W}^r_pH_\alpha[\phi] &=\left\{ f\in \widetilde{W}^r_p:
\omega_\alpha \left(\widetilde{f}^{(r)},
\delta \right)_p = O[\phi(\delta)] \right\}, \\
{W}^r _pE[\xi] &=\Bigl\{ f\in W^r_p: E_n \left(f^{(r)}\right)_p =
O[\xi(1/n)] \Bigr\}.
\end{align*}
We will study more general classes  such that  $W^r_p$,
$\widetilde{W}^r_p$, ${W}^r_pH_\alpha[\varphi]$,
$\widetilde{W}^r_pH_\alpha[\phi]$, ${W}^r E_p[\xi]$ are their particular cases.
\\
{\sf Transformed Fourier series}. Let $L_p=L_p[0, 2\pi] \,\,(1\le
p< \infty)$ be a space of
$2\pi$-periodic  measurable functions such that $|f|^p$ is integrable, and $L_\infty\equiv C[0,
2\pi]$ be the space of $2\pi$-periodic continuous functions with the uniform norm, that is,
$\|f\|_\infty=\max\limits\left\{|f(x)|,
0\le x\le 2\pi\right\}.$
\\
Let the Fourier series of a summable function
$f(x)$ be written as
\begin{equation}\label{fourier}
f(x) \sim \sigma(f):= \frac{a_0 (f)}{2}+
\sum\limits_{\nu=1}^\infty \left( a_\nu(f) \cos \nu x + b_\nu(f)
\sin \nu x\right) \equiv \sum\limits_{\nu=0}^\infty A_\nu(f,x).
\end{equation}
{\it The transformed Fourier series}  for series (\ref{fourier}) is defined as follows:
 $$\sigma(f,\lambda,\beta):=\sum\limits_{\nu=1}^\infty
\lambda_\nu \left[ a_\nu \cos\left( \nu x
+\frac{\pi\beta}{2}\right) + b_\nu \sin \left( \nu x
+\frac{\pi\beta}{2}\right) \right],$$ where
$\beta\in \mathbf{R}$  and
 $\lambda=\left\{\lambda_n\right\}$ is a given sequence of positive numbers.

This definition is well-known in the literature (see, for example, \cite[ch. 12 \S 8-9]{zygm}, \cite{og}, \cite{potapov},
\cite{simonov}, \cite{step}, \cite{cohen}, \cite{sch}, \cite{sz} and
remark 3 in this paper). We also note that $\sigma(f,\lambda,\beta)$
coincides up to notations with the Fourier series of so called
$f^\psi_\beta$-derivatives, using the terminology of \cite[p. 132]{step}.

Studies of the  transformed Fourier series  are naturally related to the problems of Fourier
multipliers theory (see \cite{nikol'skij}, \cite[V. 1, Ch. III]{zygm}, \cite{berg}, \cite[Chapter 7]{trigub}),
summability methods (see \cite[V. 1, Ch. III]{zygm}, \cite[Chapter 1.2]{butzbook}, \cite[Chapter 8]{trigub})
and\footnote{See also references to \cite[\S 13, Ch. II]{baritext}.} the so-called fractional
Sobolev classes or the Weyl classes \cite[V. 1, Ch. III]{step}.

The function class
\begin{equation*}
W_p^{\lambda, \beta} =\Bigl\{ f\in L_p :\quad \exists \,\,g
 \in L_p, \quad \sigma(g)= \sigma(f,\lambda,\beta) \Bigr\}
\end{equation*}
is called {\it the Weyl class} (see for example \cite{potapov}, \cite{step}, \cite{ps1}).
It is named so, because for
$\lambda_n = n^r, r>0$ and
$\beta = r$ the class
  $W_p^{\lambda, \beta}$
  coincides with the class $W^r_p$, which is defined
in terms of fractional derivatives $f^{(r)}$ in the Weyl sense (\cite[V. 2, Ch.
XII]{zygm}).
  In the case of $\lambda_n = n^r, r>0$ and $\beta = r+1$ the class $W_p^{\lambda, \beta}$
coincides with the class $\widetilde{W}^r_p$.
   A function $g(x) \sim \sigma(f,\lambda,\beta)$
   is called {\it the $(\lambda,\beta)$-derivative} of a function $f(x)$ and is denoted  by
$f^{(\lambda,\beta)}(x)$.
 Using the terminology of \cite{step}, we have
$f^\psi_\beta=f^{(\lambda,\beta)}$ for $\psi^{-1}(k)=\lambda_k$.
\\
{\sf The generalized Weyl-Nikolskii class}. In the definition of this functional class we use {\it the
modulus of smoothness} concept $\omega_{\alpha} (f,\delta)_p$ of {\it fractional
\footnote{ The term
"fractional" can be found in earlier papers
 (\cite{butz} and \cite{tab}) which used this  definition.
As in the case of fractional derivatives, the positive number\, $\alpha$ that defines  the modulus
order is not necessarily rational.}
order} of a function $f(x) \in L_p$, i.e.,
$$\omega_{\alpha} (f,\delta)_p
 =\sup_{|h|\le \delta }
 \left\|
\triangle_{h}^{\alpha} f(x) \right\|_p,$$ where
 $$
\triangle_{h}^{\alpha} f(x) =\sum\limits_{\nu=0}^\infty(-1)^{\nu}
\binom{\alpha}{\nu} f(x+(\alpha-\nu)h), \qquad \alpha>0$$
is the $\alpha$-th difference\footnote{As usual,
 $\binom{\beta}{\nu} =
\frac{\beta(\beta-1)\cdots(\beta-\nu+1)}{\nu!}$ for $\nu
> 1$, $\binom{\beta}{\nu}=\beta$ for $\nu = 1$, and
 $\binom{\beta}{\nu} = 1 $ for $\nu = 0$.}
 of a function $f$ with step $h$ at the point $x$.
It is clear that for $\alpha\in \mathbf{N}$ this definition is the same as (\ref{mod}).
\\
Let $\Phi_\alpha\,\,(\alpha>0)$
 be the class of functions
 $\varphi(\delta)$, defined and non-negative on  $(0,\pi]$ such that

1. $\varphi(\delta)\to 0\,\,\,(\delta\to 0)$,

2. $\varphi(\delta)$  is non-decreasing,

3. $\delta^{-\alpha}\varphi(\delta)$ is non-increasing.
\\
For functions $\varphi\in \Phi_\alpha, \alpha>0$ and for $\lambda=\left\{\lambda_n\right\}$ we define
{\it the generalized Weyl-Nikolskii class}  similarly to the classes
 ${W}^r_pH^\alpha[\varphi]$ and
$\widetilde{W}^r_pH^\alpha[\phi]$ (see, for example, \cite{ps1}):
\begin{equation*}
W_p^{\lambda, \beta} H_\alpha[\varphi] = \Bigl\{ f\in
W_p^{\lambda, \beta} :\,\,\, \omega_\alpha\left(
f^{(\lambda,\beta)},\delta \right)_p=
O\left[\varphi(\delta)\right], \delta \to +0
 \Bigr\}.
\end{equation*}
It is clear that if  $\lambda_n = n^r, r>0$ and $\beta = r$, then \,
 $W_p^{\lambda, \beta} H_\alpha[\varphi]\equiv  {W}^r_pH_\alpha[\varphi]$;
and if $\lambda_n = n^r, r>0$ and $\beta = r+1$, then \,
 $W_p^{\lambda, \beta} H_\alpha[\varphi]\equiv
 \widetilde{W}^r_pH_\alpha[\varphi]$.

In case $\lambda_n\equiv 1$ and $\beta=0$ the class
 $W_p^{\lambda, \beta} H_\alpha[\varphi]$ coincides with  {\it the generalized Lipschitz class}
 $H_\alpha^\varphi$, i.e.,
\begin{equation*} H_\alpha^p[\varphi] = \Bigl\{ f\in L_p :
\omega_\alpha(f,\delta)_p= O\left[\varphi(\delta)\right]\quad
\delta \to +0 \Bigr\}.
\end{equation*}
In particular, for $0<\gamma\le 1$,
\begin{equation*}
\textnormal{Lip} \,(\gamma, L_p) \equiv H_1^p \,[\delta^\gamma]=
\Bigl\{ f\in L_p : \omega_1(f,\delta)_p=
O\left[\delta^\gamma\right]\quad \delta \to +0 \Bigr\}.
\end{equation*}
{\sf The problem setting and the structure of the paper}. In this paper, we obtain embedding
theorems for the Weyl classes $W_p^{\lambda, \beta}$, for the generalized Weyl-Nikolskii classes
$W_p^{\lambda, \beta} H_\alpha[\varphi]$ and for the generalized Lipschitz classes
$H_\gamma^p[\omega]$.
We show how the  parameters\, $\alpha$ and $\gamma$ are related to each other depending
on the behavior of the sequence $\{\lambda_n\}$ and on the choice of the metric $L_p$.

This paper  is organized as follows. In section 2 we formulate the main theorem.
Sections 3 and  4 contain the proofs of the sufficiency and necessity parts of the main theorem
respectively.
In section 5 we  prove several corollaries. In particular, we describe the
difference in results for metrics
 $L_p, 1<p<\infty$  and  $L_p, p=1,\infty$.
Also, the estimates of
$\omega_\gamma(f^{(r)}, \delta)_p$ and
$\omega_\gamma(\widetilde{f}^{(r)}, \delta)_p$ are written in terms of
$\omega_\beta(f, \delta)_p$ for different values of $r, \gamma$, and
$\beta$. The concluding remarks are given in section 6.

\section{Embedding theorems for the generalized Lipschitz and \\
 Weyl-Nikolskii   classes}

For $\lambda=\left\{\lambda_n\right\}_{n\in \mathbf{N}}$ we define
$\triangle\lambda_n:=\lambda_n- \lambda_{n+1}$;
$\triangle^2\lambda_n:=\triangle\left(\triangle\lambda_n\right)$.

\begin{thr} \label{osnth}
Let\, $\theta=\min (2,p)$,  $\alpha\in  \mathbf{R}_+,$
$\beta\in \mathbf{R}$, $\rho \in  \mathbf{R}_+\cup \{0\}$
and $\lambda=\left\{\lambda_n\right\}$
be a non-decreasing sequence of
positive numbers
such that $\left\{n^{-\rho}\lambda_n\right\}$ is non-increasing.
\\
{\bf I. }\,\, If $1<p<\infty$, then
\begin{eqnarray}
\label{1} H_{\alpha+\rho}^p[\omega] \subset W_p^{ \lambda,\beta}
\quad &\Longleftrightarrow& \,\,\, \sum\limits_{n=1}^\infty
\left(\lambda_{n+1}^\theta-\lambda_{n}^\theta\right)
\omega^\theta\left(\frac{1}{n}\right)\, < \,\infty,
\\
H_{\alpha+\rho}^p[\omega] \subset W_p^{ \lambda,\beta}
H_\alpha[\varphi] \,\, &\Longleftrightarrow& \biggl\{
n^{-\alpha\theta} \sum\limits_{\nu=1}^n \nu^{(\rho+\alpha)\theta}
\left(\nu^{-\rho\theta}\lambda_{\nu}^\theta-
(\nu+1)^{-\rho\theta}\lambda_{\nu+1}^\theta \right)
\omega^\theta\left(\frac{1}{\nu}\right) \nonumber
\\
&+& \sum\limits_{\nu=n+2}^\infty
\left(\lambda_{\nu+1}^\theta-\lambda_{\nu}^\theta\right)
\omega^\theta\left(\frac{1}{\nu}\right) + \lambda_{n+1}^\theta
\omega^\theta\left(\frac{1}{n+1}\right) \biggr\}^\frac{1}{\theta}
\nonumber
\\
&=& \label{2} \quad  O\left[\varphi\left(\frac{1}{n+1}
\right)\right],
\\
\label{3} W_p^{\lambda,\beta}
 \subset H_{\alpha+\rho}^p[\omega]
\quad &\Longleftrightarrow& \quad \frac{1}{\lambda_{n}}= O
\left[\omega\left(\frac{1}{n} \right)\right],
\\
\label{4} W_p^{ \lambda,\beta} H_\alpha[\varphi]
 \subset H_{\alpha+\rho}^p[\omega]
\quad &\Longleftrightarrow& \quad \frac{\varphi\left(\frac{1}{n}
\right)}{\lambda_{n}}= O \left[\omega\left(\frac{1}{n}
\right)\right].
\end{eqnarray}
{\bf II. } \,\,Let $p=1$ or $p=\infty$.
\\
\textnormal{(a)} If   $\triangle^2\lambda_n\ge 0$ or
$\triangle^2\lambda_n\le 0$, then
\begin{eqnarray}\label{5}
H_{\alpha+\rho}^p[\omega] \subset W_p^{ \lambda,\beta}
&\Longleftrightarrow& | \cos \frac{\beta\pi}{2} |
\sum\limits_{n=1}^\infty  \left(\lambda_{n+1}-\lambda_{n}\right)
\omega\left(\frac{1}{n}\right) \nonumber
\\
&+& | \sin \frac{\beta\pi}{2} | \sum\limits_{n=1}^\infty
\lambda_{n} \frac{\omega\left(\frac{1}{n}\right)}{n} \, <
\,\infty;
\end{eqnarray}
and if, additionally, for some  $\tau > 0$ the following inequality holds,
$$\triangle^2\left(\frac{\lambda_n}{n^r}\right)\ge 0 \quad\mbox{for}\quad
r=\rho+\tau\, sign \left|\sin\frac{(\beta-\rho)\pi}{2}\right|,$$
then
\begin{eqnarray}\label{6}
H_{\alpha+r}^p[\omega] \subset W_p^{ \lambda,\beta}
H_\alpha[\varphi] &\Longleftrightarrow& n^{-\alpha}
\sum\limits_{\nu=1}^n \nu^{r+\alpha} \left(\nu^{-r}\lambda_{\nu}-
(\nu+1)^{-r}\lambda_{\nu+1} \right)
\omega\left(\frac{1}{\nu}\right) \nonumber
\\
&+& | \cos \frac{\beta\pi}{2} | \sum\limits_{\nu=n+2}^\infty
\left(\lambda_{\nu+1}-\lambda_{\nu}\right)
\omega\left(\frac{1}{\nu}\right) \nonumber
\\
&+& | \sin \frac{\beta\pi}{2} | \sum\limits_{\nu=n+2}^\infty
 \lambda_\nu
\frac{\omega\left(\frac{1}{\nu}\right)}{\nu}+ \lambda_{n+1}
\omega\left(\frac{1}{n+1}\right) \nonumber
\\&=& \quad
O\left[\varphi\left(\frac{1}{n+1}\right)\right].
\end{eqnarray}
\\
\textnormal{(b)} If  for $\beta= 2k, \, k\in \mathbf{Z}$, the condition
$\triangle^2\left(1/\lambda_n\right)\ge 0$ holds, and for $\beta\ne 2k, \, k\in \mathbf{Z}$ conditions
$\triangle^2\left(1/\lambda_n\right)\ge 0$ and
$\sum\limits_{\nu=n+1}^\infty \frac{1}{\nu\lambda_\nu} \le
\frac{C}{\lambda_n}$ are fulfilled, then
\begin{equation}\label{7}
W_p^{\lambda,\beta} \subset
 H^p_{\alpha+\rho}[\omega]
\quad \Longleftrightarrow \,\, \frac{1}{\lambda_{n}}= O
\left[\omega\left(\frac{1}{n} \right)\right];
\end{equation}
and if, additionally, for some  $\tau > 0$ the following inequality holds,
$$
\triangle^2\left(\frac{n^r}{\lambda_n}\right)\ge 0
\quad\mbox{or}\quad
\triangle^2\left(\frac{n^r}{\lambda_n}\right)\le 0
 \quad\mbox{for}\quad
r=\rho+\tau\, sign \left|\sin\frac{(\beta-\rho)\pi}{2}\right|,$$
then
\begin{equation}\label{8} W_p^{ \lambda,\beta} H_\alpha[\varphi]
 \subset H_{\alpha+r}^p[\omega]
\quad \Longleftrightarrow \quad \frac{\varphi\left(\frac{1}{n}
\right)}{\lambda_{n}}= O \left[\omega\left(\frac{1}{n}
\right)\right].
\end{equation}
\end{thr}

\section{Proof of sufficiency in Theorem 1.}

We will use the following notations.
\\
Let series (\ref{fourier}) be the Fourier series of a function $f(x)\in L$.
Then $S_n(f)$ denotes the $n$-th partial sum of series (\ref{fourier}), $V_n(f)$ denotes the de la Vall{\'e}e-Poussin sum and $K_n(x)$ is the Fej{\'e}r kernel, i.e.,
$$ \quad S_n(f)=\sum\limits_{\nu=0}^n A_\nu(x), \quad V_n(f)=\frac{1}{n}\sum\limits_{\nu=n}^{2n-1}
S_\nu(f), \quad K_n(x) = \frac{1}{n+1} \sum\limits_{\nu=0}^n
\left( \frac{1}{2}+\sum\limits_{m=1}^\nu \cos m x \right).$$
The following lemmas will play an important role in the proof of the main theorem.

\begin{lemma}  \label{my}
If  $f(x)\in L_p,$ $1\le p\le \infty$ and $\alpha>0$, then
\begin{equation} \label{m1}
C_1(p,\alpha) \omega_\alpha\Bigl(f,\frac{1}{n}\Bigr)_p \le \left(
n^{-\alpha}\left\| V_n^{(\alpha)}(f,x)) \right\|_p +
\left\|f(x)-V_n(f,x)\right\|_p \right) \le C_2(p,\alpha)
\omega_\alpha\Bigl(f,\frac{1}{n}\Bigr)_p.
\end{equation}
If $f(x)\in L_p,$ $1< p < \infty$, then
\begin{equation} \label{m2}
C_1(p,\alpha) \omega_\alpha\Bigl(f,\frac{1}{n}\Bigr)_p \le \left(
n^{-\alpha}\left\| S_n^{(\alpha)}(f,x)) \right\|_p +
\left\|f(x)-S_n(f,x)\right\|_p \right) \le C_2(p,\alpha)
\omega_\alpha\Bigl(f,\frac{1}{n}\Bigr)_p.
\end{equation}
\end{lemma}
 {\em Proof of Lemma  \ref{my}}. The estimate of
$\omega_\alpha\Bigl(f,\frac{1}{n}\Bigr)_p$ from above follows from the inequality \,
(see \cite{butz})
$\omega_\alpha\Bigl(T_n,\frac{1}{n}\Bigr)_p \le C(p,\alpha)
n^{-\alpha}\left\| T_n^{(\alpha)} \right\|_p$, where $T_n$ is a trigonometric polynomial of order $n$. Indeed,
\begin{eqnarray*}
\omega_\alpha\Bigl(f,\frac{1}{n}\Bigr)_p &\le&  C(p,\alpha) \left(
\omega_\alpha\Bigl(T_n,\frac{1}{n}\Bigr)_p +\left\| f-T_n
\right\|_p \right)
\\
&\le&  C(p,\alpha) \left( n^{-\alpha}\left\| T_n^{(\alpha)}
\right\|_p + \left\| f-T_n \right\|_p  \right).
\end{eqnarray*}
To estimate $\omega_\alpha\Bigl(f,\frac{1}{n}\Bigr)_p$ from below,
we will use the generalized Nikol'skii-Stechkin inequality
(see  \cite{tab})
 $n^{-\alpha}\left\| T_n^{(\alpha)} \right\|_p
\le  C(p,\alpha) \omega_\alpha\Bigl(T_n,\frac{1}{n}\Bigr)_p$
and the generalized Jackson inequality (see \cite{butz} for $\alpha>0$)
\begin{equation}\label{11}
 E_n(f)_p\le  C(\alpha) \,\,
\omega_{\alpha} \left(f,\frac{1}{n+1}\right)_p.
\end{equation}
Also, it is well known that the Vall{\'e}e-Poussin mean is the near best approximant, i.e.,
\begin{equation}\label{20}
\left\|f- V_n(f)\right\|_p\le  C E_n(f)_p.
\end{equation}
Then
\begin{eqnarray*}
n^{-\alpha}\left\| V_n^{(\alpha)}(f,x) \right\|_p +
\left\|f(x)-V_n(f,x)\right\|_p &\le& C(p,\alpha) \left(
\omega_\alpha\Bigl(V_n,\frac{1}{n}\Bigr)_p + E_n(f)_p \right)
\\
&\le& C(p,\alpha) \left( \omega_\alpha\Bigl(f,\frac{1}{n}\Bigr)_p
+ \omega_\alpha\Bigl(f-V_n,\frac{1}{n}\Bigr)_p \right)\\ &\le&
C(p,\alpha) \omega_\alpha\Bigl(f,\frac{1}{n}\Bigr)_p,
\end{eqnarray*}
i.e., (\ref{m1}) is proved.
Using
\begin{equation}\label{zv}
\left\|f- S_n(f)\right\|_p\le  C(p)\, E_n(f)_p
\end{equation}
 for $1< p < \infty$, we obtain (\ref{m2}) analogously.
 Lemma \ref{my} is proved.

We note that (\ref{m1}) and (\ref{m2}) are the realization results for modulus of smoothness
(see the original paper \cite{ditzian} by Z. Ditzian, V.  Hristov, K. Ivanov).

\begin{lemma}   \textnormal{(\cite{stechkin})}.
\label{4l22} Let $f(x)\in L_p,$ $p=1,\infty$ and let  the condition
$\sum\limits_{n=1}^{\infty}n^{-1}E_n(f)_p<\infty$ hold. Then
$\tilde{f}(x)\in L_p$ and
\begin{equation*}
E_n(\tilde{f})_p \le C \left( E_n(f)_p +
\sum\limits_{k=n+1}^{\infty} k^{-1} E_k(f)_p  \right), \qquad
n\in  \mathbf{N}.
\end{equation*}
\end{lemma}

\begin{lemma} \label{my2}
Let $p=1,\infty$ and let $\{\lambda_n\}$ be monotone  concave (or convex) sequence. Let
\begin{eqnarray*}T_n(x)& = &\sum\limits_{\nu=0}^{n} a_\nu \cos \nu x +
b_\nu \sin \nu x, \\
T_n(\lambda, x)& = &\sum\limits_{\nu=0}^{n} \lambda_\nu\left(
a_\nu \cos \nu x + b_\nu \sin \nu x\right).
\end{eqnarray*}
Then for any integer $M>N+2$ we have
$$\left\| T_M(\lambda, x) - T_N(\lambda, x) \right\|_p \le \mu(M,N)
\left\| T_M(x) - T_N(x) \right\|_p,$$ where 
$$ \mu(M,N) =
\begin{cases}
2M(\lambda_M- \lambda_{M-1})+ \lambda_{N+1} -
(N+1)(\lambda_{N+2}- \lambda_{N+1}),
 &\text{if $\lambda_n \uparrow$
\textnormal{(}$n\uparrow$\textnormal{)}, $\triangle^2\lambda_n\ge 0$};\\
2\lambda_M+(N+1)(\lambda_{N+2}-\lambda_{N+1})-\lambda_{N+1},
 &\text{if $\lambda_n \uparrow$
\textnormal{(}$n\uparrow$\textnormal{)}, $\triangle^2\lambda_n\le 0$};\\
(N+1)(\lambda_{N+1}-\lambda_{N+2})+\lambda_{N+1},
 &\text{if $\lambda_n \downarrow$
\textnormal{(}$n\uparrow$\textnormal{)}, $\triangle^2\lambda_n\ge
0$},
 \end{cases}
$$
for $M=N+1$, $ \mu(M,N) = \lambda_M$, and for $M=N+2$, $ \mu(M,N) =
2|\lambda_{N+1}-\lambda_{N+2}|+\lambda_{N+2}.$
\end{lemma}
{\em Proof of Lemma \ref{my2}}. First we consider the case when $M>N+2$.
 Applying twice Abel's transformation, we write
\begin{eqnarray*}
&\,& \left\| T_M(\lambda, x) - T_N(\lambda, x) \right\|_p =
\frac{1}{\pi} \Bigl\| \int\limits_{-\pi}^\pi \left(T_M -
T_N\right) (x+u) \sum\limits_{\nu=N+1}^{M} \lambda_\nu \cos \nu u
\, du  \Bigr\|_p
\\
&=&  \Bigl\|\frac{1}{\pi} \int\limits_{-\pi}^\pi \left(T_M -
T_N\right) (x+u) \Bigl\{ \sum\limits_{\nu=N+1}^{M-2}
\left(\lambda_{\nu}- 2 \lambda_{\nu+1} + \lambda_{\nu+2}\right)
(\nu+1) K_\nu(u)
\\
&+& \left(\lambda_{N+2}- \lambda_{N+1}\right) (N+1) K_N(u) +
\left(\lambda_{M-1}- \lambda_{M}\right) M K_{M-1}(u)\Bigr\}du +
\lambda_{M}
 \left(T_{M}-T_N\right) \Bigr\|_p
\\  &\le&
\left\| T_M(x) - T_N(x) \right\|_p \Bigl\{
\sum\limits_{\nu=N+1}^{M-2}\left|\lambda_{\nu}- 2 \lambda_{\nu+1}
+ \lambda_{\nu+2}\right| (\nu+1) + \left|\lambda_{M-1}-
\lambda_{M}\right| M + \lambda_{M} \Bigr\}
\\
&=:&
 \left\| T_M(x) - T_N(x) \right\|_p I(M,N).
\end{eqnarray*}
We estimate $I(M,N)$  in the case of  $\lambda_n \uparrow$ ($n\uparrow$),
$\triangle^2\lambda_n\ge 0$. Then
\begin{eqnarray*}
I(M,N) &=& \sum\limits_{\nu=N+1}^{M-2} \left(\lambda_{\nu}- 2
\lambda_{\nu+1} + \lambda_{\nu+2}\right) (\nu+1) +
\left(\lambda_{M}- \lambda_{M-1}\right) M + \lambda_{M}
\\&=&
 -(N+1)\left(\lambda_{N+2}- \lambda_{N+1}\right) +
\left(\lambda_{N+1}- \lambda_{M-1}\right) +\lambda_{M}+
(2M-1)\left(\lambda_{M}- \lambda_{M-1}\right)
\\
&=&
 -(N+1)\left(\lambda_{N+2}- \lambda_{N+1}\right) +
\lambda_{N+1}+2M\left(\lambda_{M}- \lambda_{M-1}\right).
\end{eqnarray*}
If $\lambda_n \uparrow$ ($n\uparrow$), $\triangle^2\lambda_n\le
0$, then
$$
I(M,N) = -(N+1)\left(\lambda_{N+1}- \lambda_{N+2}\right) + \left(
\lambda_{M-1}-\lambda_{N+1}\right) +\lambda_{M}+ \left(\lambda_{M}-
\lambda_{M-1}\right).
$$
Finally, if $\lambda_n \downarrow$ ($n\uparrow$),
$\triangle^2\lambda_n\ge 0$, then
$$
I(M,N) = -(N+1)(\lambda_{N+2}-\lambda_{N+1})+\lambda_{N+1}.
$$

The estimate for the case of  $M=N+1$ is trivial and for the case of $M=N+2$
immediately follows from the equation
\begin{eqnarray*} T_M(\lambda,x) - T_N(\lambda,x)&=&
\frac{\lambda_{M}}{\pi}
\int_{-\pi}^\pi (T_M - T_N )(x+u) \Bigl(\sum_{m=N+1}^{N+2}  \cos mu\Bigr) du\\
&+&  \frac{\lambda_{N+1} - \lambda_{N+2}}{\pi} \int_{-\pi}^\pi (T_M
- T_N )(x+u) \cos (N+1)u du.
\end{eqnarray*}
The proof of  Lemma \ref{my2} is complete.

\begin{lemma} \label{my222}
Let $p=1,\infty$. If $ \,\,T_{2^n, \,2^{n+1}}(x) =
\sum\limits_{\nu=2^n}^{2^{n+1}} \Bigl( c_\nu \cos \nu x + d_\nu
\sin \nu x \Bigr),$ then
\begin{equation}\label{trig}
C_1 \left\| \widetilde{T}_{2^n, \,2^{n+1}}(x) \right\|_p \le
\left\| T_{2^n, \,2^{n+1}}(x) \right\|_p \le C_2  \left\|
\widetilde{T}_{2^n, \,2^{n+1}}(x) \right\|_p.
\end{equation}
\end{lemma}
{\em Proof of Lemma \ref{my222}}.
 We rewrite $T_{2^n, \,2^{n+1}}(x)$ in the following way
$$
T_{2^n, \,2^{n+1}}(x) =
\sum\limits_{\nu=2^n}^{2^{n+1}}\frac{1}{\nu} \Bigl( \nu c_\nu
\cos \nu x + \nu d_\nu \sin \nu x \Bigr).
$$
Applying Lemma \ref{my2} and the Bernstein inequality, we have
\begin{eqnarray*}
\left\| T_{2^n, \,2^{n+1}}(x) \right\|_p &\le& C \frac{1}{2^{n}}
\left\| \sum\limits_{\nu=2^n}^{2^{n+1}} \Bigl( \nu  c_\nu \cos
\nu x + \nu d_\nu \sin \nu x \Bigr)
\right\|_p \\
&=& C  \frac{1}{2^{n}} \left\| \left(
\sum\limits_{\nu=2^n}^{2^{n+1}}
 -d_\nu \cos \nu x + c_\nu \sin \nu x \right)'
\right\|_p \\&\le& C \left\| \widetilde{T}_{2^n, \,2^{n+1}}(x)
\right\|_p.
\end{eqnarray*}
Similar reasoning for $\widetilde{T}_{2^n, \,2^{n+1}}(x)$ allows us to  obtain the left-hand side inequality in
(\ref{trig}). The proof of lemma \ref{my222} is now complete.

\begin{center}
\bf{Sufficiency in (\ref{1}) - (\ref{8}).}
\end{center}

\textbf{I.} $1<p<\infty$. In this case, for $\lambda_n\equiv 1$, the Riesz inequality (\cite[V. 1,
p. 253]{zygm}) $\|\tilde{f}\|_p\le C(p) \|f\|_p$ implies
\begin{equation}\label{9}
\|f^{(\lambda,\beta)}\|_p \le C(p,\beta)\|f\|_p
\end{equation}
(Here and  henceforth, by $C(s,t,\cdots)$ we understand positive constants that depend only on $s,t,\cdots$ and
in general, may be different in different inequalities).

Let the series in the  right part of (\ref{1}) be convergent and $f\in H^p_{\alpha+\rho}[\omega].$
We will use the following representation
$$
\lambda_{2^n}^\theta =
\left\{
  \begin{array}{ll}
    \lambda_{1}^\theta + \sum\limits_{\nu=2}^{n+1} \left(
\lambda_{2^{\nu-1}}^\theta - \lambda_{2^{\nu-2}}^\theta \right), \qquad& \hbox{if}\qquad n\ge 1; \\
    \lambda_{1}^\theta, & \hbox{if}\qquad {n=0.}
  \end{array}
\right.
$$
Applying Minkowski's inequality, we get $\Bigl($ here and further
 $\triangle_1:=A_1(f,x), \triangle_{n+2}:=\sum\limits_{\nu=2^n+1}^{2^{n+1}} A_\nu(f,x)$, where
$A_\nu(f,x)$ is from (\ref{fourier}) $ \Bigr)$
\begin{eqnarray} \label{AA}
I_1  &:=&   \left\{\int\limits_0^{2\pi}\left[
\sum\limits_{n=1}^{\infty}
 \lambda_{2^{n-1}}^2 \triangle^2_n\right]^\frac{p}{2}\,dx
\right\}^\frac{\theta}{p} \nonumber
\\
&\le& C(p)\left(
 \lambda_{1}^\theta
\left\{ \int\limits_0^{2\pi}\left[ \sum\limits_{n=1}^{\infty}
\triangle^2_n\right]^\frac{p}{2}\,dx \right\}^\frac{\theta}{p} +
\sum\limits_{s=2}^{\infty} \left( \lambda_{2^{s-1}}^\theta
-\lambda_{2^{s-2}}^\theta\right) \left\{
\int\limits_0^{2\pi}\left[ \sum\limits_{n=s}^{\infty}
\triangle^2_n\right]^\frac{p}{2}\,dx \right\}^\frac{\theta}{p}
\right)^\frac{1}{\theta}.
\end{eqnarray}
By the Littlewood-Paley theorem  (see, for example, \cite{zygm}, p. 349) and inequality (\ref{zv}),
we obtain
\begin{equation}\label{zv1}
I_1\le C(p)\left\{
 \lambda_{1}^\theta
 \|f\|^\theta_p +
\sum\limits_{s=1}^{\infty} \left( \lambda_{2^{s}}^\theta
-\lambda_{2^{s-1}}^\theta\right) E_{2^{s-1}}^\theta(f)_p
\right\}^\frac{1}{\theta}.
\end{equation}
Then, both the generalized Jackson inequality (\ref{11}) and  the condition
 $f\in H^p_{\alpha+\rho}[\omega]$
imply $I_1<\infty.$ Thus, there exists a function  $g\in L_p$ with the Fourier series
\begin{equation}\label{12}
\sum\limits_{n=1}^{\infty} \lambda_{2^{n-1}} \triangle_n,
\end{equation}
and $\|g\|_p\le C(p) I_1$.
We write series  (\ref{12}) in the form of
 $\sum\limits_{n=1}^{\infty}
\gamma_{n} A_n(f,x)$, where $\gamma_{i}:=\lambda_i, \,i=1,2$ and
$\gamma_{\nu}:= \lambda_{2^{n}}$ for $2^{n-1}+1\le \nu \le 2^{n} \,(n=2,3,\cdots).$
 Further, we consider the series
\begin{equation}\label{13}
\sum\limits_{n=1}^{\infty}   \lambda_{n}A_n(f,x) =
\sum\limits_{n=1}^{\infty}   \gamma_n \Lambda_{n}A_n(f,x),
\end{equation}
where $\Lambda_{1}=\Lambda_{2}=1$, $\Lambda_{\nu} :=
\lambda_{\nu}/\gamma_{n} = \lambda_{\nu}/\lambda_{2^{n}}$ for
$2^{n-1}+1\le \nu \le 2^{n} \,(n=2,3,\cdots).$

Since the sequence $\left\{\Lambda_{n}\right\}$ satisfies the conditions of the Marcinkiewicz
multiplier theorem (\cite{zygm}, p.346), series (\ref{13}) is the Fourier series of a
function $f^{(\lambda, 0)} \in L_p$ and $\|f^{(\lambda, 0)}\|_p \le C(p) \|g\|_p.$
 Then from inequalities (\ref{11}), (\ref{9}) and  (\ref{zv1}) we get
\begin{eqnarray}\label{14}
\|f^{(\lambda,\beta)}\|_p &\le& C(p,\beta) \left\{
 \lambda_{1}^\theta
 \|f\|^\theta_p +
\sum\limits_{s=1}^{\infty} E_{2^{s-1}}^\theta(f)_p
\sum\limits_{n=2^{s-1}}^{2^{s}-1} \left( \lambda_{n+1}^\theta
-\lambda_n^\theta\right)
  \right\}^\frac{1}{\theta}
\\
&\le& C(p, \beta, \alpha, \rho) \left\{
 \lambda_{1}^\theta
 \|f\|^\theta_p +
\sum\limits_{n=1}^{\infty} \left( \lambda_{n+1}^\theta
-\lambda_n^\theta\right)
 \omega_{\alpha+\rho}^\theta \left(f,\frac{1}{n}\right)_p
 \right\}^\frac{1}{\theta},
\nonumber
\end{eqnarray}
i.e., the sufficiency in (\ref{1}) is proved.

Let the relation in the right-hand side in (\ref{2}) hold, and $f\in
H^p_{\alpha+\rho}[\omega]$.
Let us prove $f\in W^{\lambda,\beta}_pH_{\alpha}[\varphi].$ First, we estimate
$\omega_\alpha(f^{(\lambda,\beta)},\frac{1}{n})_p.$
By Lemma \ref{my},
\begin{equation}\label{16}
\omega_\alpha\left(f^{(\lambda,\beta)},\frac{1}{n}\right)_p \le
C(p,\alpha) \Bigl( \|f^{(\lambda,\beta)} -
S_n(f^{(\lambda,\beta)}) \|_p+ n^{-\alpha}
\|S^{(\alpha)}_n(f^{(\lambda,\beta)}) \|_p
 \Bigr).
\end{equation}
Using (\ref{14}) for the function $\left(f - S_n\right)$, we have ($[a]$ is the integer part of
$a$)
\begin{align} \label{17}
\|f^{(\lambda,\beta)} - S_n(f^{(\lambda,\beta)})\|_p \le C(p,
\beta, \alpha, \rho) \Biggl\{
 \lambda_{n+1}^\theta
\left\|f - S_n\right\|^\theta_p \, &+ \,\,
E_{\left[\frac{n}{2}\right]}^\theta(f)_p \sum\limits_{s=1}^{2n}
\left( \lambda_{s+1}^\theta -\lambda_s^\theta\right) \nonumber
\\
&+        \,\, \sum\limits_{s=n+1}^{\infty} \left(
\lambda_{s+1}^\theta -\lambda_s^\theta\right)
E_{\left[\frac{s}{2}\right]}^\theta(f)_p
  \Biggr\}^\frac{1}{\theta}
\nonumber
\\
\le              \, C(p, \beta, \alpha, \rho) \Biggl\{
 \lambda_{n+1}^\theta
\omega^\theta_{\alpha+\rho} \left(f,\frac{1}{n}\right)_p \,\, &+
\, \left. \sum\limits_{\nu=n+1}^{\infty} \left(
\lambda_{\nu+1}^\theta -\lambda_\nu^\theta\right)
 \omega_{\alpha+\rho}^\theta \left(f,\frac{1}{\nu}\right)_p
 \right\}^\frac{1}{\theta}.
\end{align}

Further, we estimate the second term of (\ref{16}). Let $m$ be an integer such that
$2^m\le n+1< 2^{m+1}.$ We will use the identity
$$2^{-s\rho\theta}\lambda_{2^s}^\theta =
2^{-(m+1)\rho\theta} \lambda_{2^{m+1}}^\theta +
\sum\limits_{\nu=s}^{m} \left(
2^{-\nu\rho\theta}\lambda_{2^\nu}^\theta -
2^{-(\nu+1)\rho\theta}\lambda_{2^{\nu+1}}^\theta \right).$$
Then using Lemmas \ref{my} and \ref{my2},  we follow the proof of typical estimates
(\ref{AA})-(\ref{14}). Then we get
\begin{eqnarray} \label{18}
n^{-\alpha} \|S^{(\alpha)}_n\left(f^{(\lambda,\beta)}\right) \|_p
&\le&
 C(p, \beta, \alpha, \rho) \Biggl\{
  \lambda_{n+1}^\theta
\omega^\theta_{\alpha+\rho}(f,\frac{1}{n})_p \nonumber
\\ &+&
n^{-\alpha\theta} \sum\limits_{\nu=1}^{n} \left(
\nu^{-\rho\theta}\lambda_{\nu}^\theta -
(\nu+1)^{-\rho\theta}\lambda_{\nu+1}^\theta \right)
\nu^{(\rho+\alpha)\theta}
 \omega_{\alpha+\rho}^\theta \left(f,\frac{1}{\nu}\right)_p
 \Biggr\}^\frac{1}{\theta}.
\end{eqnarray}
Collecting estimates (\ref{17}), (\ref{18}) and the inequality in the right-hand side of (\ref{2}), we get
 $f\in
W^{\lambda,\beta}_pH_{\alpha}[\varphi].$

Now we prove that conditions
 $\frac{1}{\lambda_{n}} = O
\left[\omega\left(\frac{1}{n}\right)\right]$ and
$\frac{\varphi\left(\frac{1}{n} \right)}{\lambda_{n}}= O
\left[\omega\left(\frac{1}{n} \right)\right]$ are sufficient for embeddings $W^{\lambda,\beta}_p \subset H^p_{\alpha+\rho}[\omega]$ and
$W^{\lambda,\beta}_p H_{\alpha}[\varphi]
 \subset H^p_{\alpha+\rho}[\omega]$, respectively.

 Using the
Littlewood-Paley and the Marcinkiewicz multiplier theorems
and the properties of the sequence $\left\{\lambda_n\right\}$ (following the proof of  (\ref{AA})-(\ref{14})), we get
\begin{eqnarray*}
 \omega_{\alpha+\rho}\left(f,\frac{1}{n}\right)_p
&\le& C(p, \alpha, \rho) \Bigl( \|f - S_{n}(f) \|_p +
n^{-(\alpha+\rho)} \|S^{(\alpha+\rho)}_{n}(f) \|_p \Bigr)
\\
&\le& C(p, \beta, \alpha, \rho) \Bigl( \lambda^{-1}_{n}
\|f^{(\lambda,\beta)} - S_{n}\left(f^{(\lambda,\beta)}\right) \|_p
+ \lambda^{-1}_{n} n^{-\alpha}
 \|S^{(\alpha)}_{n}(f^{(\lambda,\beta)}) \|_p
\Bigr).
\end{eqnarray*}
Then, by Lemma \ref{my}, we have the following inequalities
$$
\omega_{\alpha+\rho}\left(f,\frac{1}{n}\right)_p \le C(p, \beta,
\alpha, \rho) \lambda^{-1}_{n}
\omega_{\alpha}\left(f^{(\lambda,\beta)},\frac{1}{n}\right)_p \le
C(p, \beta, \alpha, \rho) \lambda^{-1}_{n}
\|f^{(\lambda,\beta)}\|_p,
$$
where the first one implies  sufficiency in (\ref{4}) and the second one shows sufficiency in (\ref{3}).

 {\bf II.} $p=1$ or $p=\infty$.
 Let the series in the right-hand side of (\ref{5}) converge, and let $f\in
H^p_{\alpha+\rho}[\omega]$. Consider the series
\begin{eqnarray} \label{19}
\cos \frac{\pi\beta}{2} V_1(\lambda,f)& -& \sin
\frac{\pi\beta}{2} \widetilde{V_1}(\lambda,f)
 \\
&+& \sum\limits_{n=1}^{\infty} \left\{ \cos \frac{\pi\beta}{2}
\left(
 V_{2^n}(\lambda,f) -  V_{2^{n-1}}(\lambda,f) \right)
-  \sin \frac{\pi\beta}{2} \left( \widetilde{V_{2^n}}(\lambda,f) -
\widetilde{V_{2^{n-1}}}(\lambda,f) \right) \right\}, \nonumber
\end{eqnarray}
where $V_{1}(\lambda,f) := \lambda_1 A_1(f,x),$
$$V_{n}(\lambda,f)
:= \sigma(\lambda, V_n(f)) =
 \sum\limits_{m=1}^{n}
 \lambda_m A_m(f,x)
+ \sum\limits_{m=n+1}^{2 n-1}
\lambda_m\left(1-\frac{m-n}{n}\right)
 A_m(f,x) \,(n\ge 2).$$
Let $M>N>0$.

Using the inequality
 $\left\|f-V_n(f)\right\|_p\le C E_n(f)_p$, the Jackson theorem (\ref{11}), and the properties of  $\left\{\lambda_n\right\}$, and following the proof of Lemma \ref{my2}, we get
\begin{eqnarray} \label{21}
A &:=&\left\|
 \sum\limits_{n=N}^{M}\left[ \cos \frac{\pi\beta}{2}
\left(
 V_{2^{n+1}}(\lambda,f) -  V_{2^{n}}(\lambda,f)
\right) - \sin \frac{\pi\beta}{2} \left(
\widetilde{V_{2^{n+1}}}(\lambda,f) -
\widetilde{V_{2^{n}}}(\lambda,f) \right) \right] \right\|_p
\nonumber \\
&\le& \sum\limits_{n=N}^{M}\left[ | \cos \frac{\pi\beta}{2}|
\left\|  V_{2^{n+1}}(f) - V_{2^{n}}(f)  \right\|_p \left(
\sum\limits_{m=2^{n}}^{2^{n+2}-1} |\triangle^2 \lambda_{m}| (m+1)
+ 2^{n+2} |\triangle \lambda_{2^{n+2}}| \right) \right.
\nonumber \\
&+& \left. | \sin \frac{\pi\beta}{2}| \left\|
\widetilde{V_{2^{n+1}}}(f) - \widetilde{V_{2^{n}}}(f)
 \right\|_p
\left( \sum\limits_{m=2^{n}}^{2^{n+2}-1} |\triangle^2
\lambda_{m}| (m+1) + 2^{n+2} |\triangle \lambda_{2^{n+2}-1}|
\right) \right]
\nonumber \\
&+& |\cos \frac{\pi\beta}{2}| \left\| \sum\limits_{n=N}^{M}
 \lambda_{2^{n+2}}
\left(  V_{2^{n+1}} -  V_{2^n}\right)(f)  \right\|_p + |\sin
\frac{\pi\beta}{2}| \left\| \sum\limits_{n=N}^{M}
 \lambda_{2^{n+2}}
\left(  \widetilde{V_{2^{n+1}}} - \widetilde{ V_{2^n}}
\right)(f)  \right\|_p
\nonumber \\
&\le& C \Biggl\{
 \lambda_{2^{N}} \left(
|\cos \frac{\pi\beta}{2}| E_{2^{N}} (f)_p + |\sin
\frac{\pi\beta}{2}| E_{2^{N}} (\tilde{f})_p \right)
\nonumber \\
&+& \sum\limits_{n=2^{N}-1}^{\infty} \left( \lambda_{n+1} -
\lambda_{n}\right) \left( |\cos \frac{\pi\beta}{2}|\,\,
\omega_{\alpha+\rho}\left(f, \frac{1}{n}\right)_p + |\sin
\frac{\pi\beta}{2}|\,\, \omega_{\alpha+\rho}\left(\tilde{f},
\frac{1}{n}\right)_p \right) \Biggr\}.
\end{eqnarray}
To complete the proof of the sufficiency part in (\ref{5}), we apply Lemma \ref{4l22}, inequality
(\ref{a3}) (see \cite{tab} for the case $k>0$),  and inequality  (\ref{11}).

 Then the convergence of the series in the right-hand side of (\ref{5}) and the condition $f\in H^p_{\alpha+\rho}[\omega]$ imply  the fact that
 the sequence
$$\left\{
 V_{2^{n}}(\lambda,\beta,f):= \cos \frac{\pi\beta}{2}
 V_{2^{n}}(\lambda,f) -  \sin \frac{\pi\beta}{2}  \widetilde{V_{2^n}}(\lambda,f)
\right\}$$ is fundamental in $L_p$. If $p=1$, since $L_1$ is complete, there exists a subsequence  $\{n_k\}$ such that $V_{2^{n_k}}(\lambda,\beta,f)$ converges almost everywhere to a function $\varphi\in L_1$.
Then from the mean convergence we obtain that, say for cosine coefficients,
$$
a_n(\varphi) =  \frac{1}{\pi} \int\limits_{-\pi}^{\pi} \varphi(x)
\cos nx \, dx = \lim\limits_{k\to\infty}\frac{1}{\pi}
\int\limits_{-\pi}^{\pi} V_{2^{n_k}}(\lambda,\beta,f) \cos nx \, dx=
a_n(f^{(\lambda,\beta)}).
$$
Therefore,
$\sigma(\varphi) = \sigma(f^{(\lambda,\beta)}).$

For $p=\infty$ the proof is similar. This completes the proof of the sufficiency part of (\ref{5}).

Let now the condition in the right-hand side of (\ref{6}) hold and
 $f\in H^p_{\alpha+r}[\omega]$. Let us estimate $\omega_\alpha\left(f^{(\lambda,\beta)},\frac{1}{n}\right)_p$ from above. By Lemma \ref{my},
$$
\omega_\alpha\left(f^{(\lambda,\beta)},\frac{1}{n}\right)_p \le
C(\alpha)\left( \left\|f^{(\lambda,\beta)}
-V_n(f^{(\lambda,\beta)})\right\|_p+ n^{-\alpha}
\left\|V_n^{(\alpha)}(f^{(\lambda,\beta)})\right\|_p \right).$$
Let us show that
\begin{eqnarray}  \label{23}
\left\|f^{(\lambda,\beta)} -V_n(f^{(\lambda,\beta)})\right\|_p
&\le& C(\beta, \alpha, r)
 \left(  \lambda_{n}
\omega_{\alpha+r}\left( f, \frac{1}{n}\right)_p +|\cos
\frac{\pi\beta}{2}|\, \sum\limits_{\nu=n+1}^{\infty} \left(
\lambda_{\nu+1} -  \lambda_{\nu}\right) \omega_{\alpha+r}\left( f,
\frac{1}{\nu}\right)_p \right.\nonumber
\\
&+& \left. |\sin \frac{\pi\beta}{2}|\,
\sum\limits_{\nu=n+1}^{\infty} \frac{\lambda_{\nu}}{\nu}
\omega_{\alpha+r}\left( f, \frac{1}{\nu}\right)_p
 \right).
\end{eqnarray}
It has been proved above that
 $$A= \left\|
  \sum\limits_{n=N}^{M}
  \left(
 V_{2^{n+1}}(f^{(\lambda,\beta)}) -  V_{2^{n}}(f^{(\lambda,\beta)})
 \right)
\right\|_p.$$
Then for  $2^m\le n<2^{m+1}$
\begin{eqnarray*}
\left\|f^{(\lambda,\beta)} -V_n(f^{(\lambda,\beta)})\right\|_p
&\le& \left\|
 V_{n}(f^{(\lambda,\beta)}) - V_{2^{m+1}}(f^{(\lambda,\beta)})
\right\|_p + \left\| \sum\limits_{\nu=m+1}^{\infty} \left(
 V_{2^\nu}(f^{(\lambda,\beta)}) - V_{2^{\nu+1}}(f^{(\lambda,\beta)})\right) \right\|_p
\\&=:& I_1+I_2.
\end{eqnarray*}
By Lemma \ref{my222}, we have
\begin{eqnarray*} I_1 &=& \left\|
\cos \frac{\pi\beta}{2} \left( V_{n}(\lambda, f) - V_{2^{m+1}}
(\lambda, f) \right) - \sin \frac{\pi\beta}{2} \left(
\widetilde{V_{n}}(\lambda, f) - \widetilde{V_{2^{m+1}}} (\lambda,
f) \right) \right\|_p
\\
&\le& \left\| V_{n}(\lambda, f) - V_{2^{m+1}}(\lambda,
f)\right\|_p + \left\|\widetilde{V_{n}}(\lambda, f) -
\widetilde{V_{2^{m+1}}} (\lambda, f) \right\|_p
 \\&\le&
 C \left\|
V_{n}(\lambda, f) - V_{2^{m+1}}(\lambda, f)\right\|_p
\end{eqnarray*}
and, by Lemma \ref{my2}, we write $ I_1 \le \lambda_{n}
\omega_{\alpha+r}\left( f, \frac{1}{n}\right)_p$.

Further, we estimate
\begin{eqnarray*}
I_2 &\le& |\cos \frac{\pi\beta}{2}|
 \left\|
 \sum\limits_{\nu=m+1}^{\infty} \left( V_{2^\nu}(\lambda, f) - V_{2^{\nu+1}}
(\lambda, f) \right) \right\|_p + | \sin \frac{\pi\beta}{2}|
 \left\|
 \sum\limits_{\nu=m+1}^{\infty}
 \left( \widetilde{V_{2^\nu}}(\lambda, f) - \widetilde{V_{2^{\nu+1}}} (\lambda, f) \right)\right\|_p.
\end{eqnarray*}
As in (\ref{21}), we write
\begin{eqnarray*}
 &&\left\|
 \sum\limits_{\nu=N}^{M} \left( V_{2^{\nu+1}} - V_{2^{\nu}}\right)
(\lambda, f)  \right\|_p \le C\left( \lambda_{2^{N}}
\omega_{\alpha+r}\left( f, \frac{1}{2^{N}}\right)_p +
 \sum\limits_{\nu=2^{N}+1}^{\infty}
 \left( \lambda_{\nu+1} -  \lambda_{\nu}\right)
\omega_{\alpha+r}\left( f, \frac{1}{\nu}\right)_p
 \right)
 ;
 \\&&
 \left\|
 \sum\limits_{\nu=N}^{M} \left(
\widetilde{V_{2^\nu}}(\lambda, f) - \widetilde{V_{2^{\nu+1}}}
(\lambda, f)
  \right) \right\|_p
  \le
   \left\|
 \sum\limits_{\nu=N}^{M}
 \lambda_{2^{\nu+2}}
 \left(
\widetilde{V_{2^{\nu}}}(f) - \widetilde{V_{2^{\nu+1}}}(f)
  \right) \right\|_p
  \\
&&+
 \sum\limits_{\nu=N}^{M}
 \left\|\widetilde{V_{2^{\nu}}}(f) - \widetilde{V_{2^{\nu+1}}}(f) \right\|_p
\left( \sum\limits_{m=2^{\nu}}^{2^{\nu+2}-1} |\triangle^2
\lambda_{m}| (m+1) + 2^{\nu+2} |\triangle \lambda_{2^{\nu+2}}|
\right) =: I_{21}+I_{22}.
\end{eqnarray*}
Applying Lemma \ref{4l22},
\begin{eqnarray*}
 I_{21}
  &\le& \lambda_{2^{N+1}}
  \left\| \widetilde{V_{2^{M+1}}}(f) - \widetilde{V_{2^N}}(f)  \right\|_p
  +
 \sum\limits_{\nu=N}^{M}  \left(\lambda_{2^{\nu+2}}-\lambda_{2^{\nu+1}}\right)
\left\| \widetilde{V_{2^{\nu+1}}}(f) -
\widetilde{V_{2^{\nu}}}(f)  \right\|_p
  \\
 &\le& C\left(
 \lambda_{2^{N+1}}
 \sum\limits_{\nu=N}^{\infty} E_{2^\nu}(f)_p
  +
 \sum\limits_{\nu=N}^{\infty}  \left(\lambda_{2^{\nu+2}}-\lambda_{2^{\nu+1}}\right)
 \sum\limits_{s=\nu}^{\infty}  E_{2^s}(f)_p\right)
\\&\le& C
 \sum\limits_{\nu=2^N}^{\infty} \frac{\lambda_{\nu}}{\nu}
\omega_{\alpha+r}\left( f, \frac{1}{\nu}\right)_p,
 \\
 I_{22}
  &\le& C
 \sum\limits_{\nu=N}^{M}  \left(\lambda_{2^{\nu+3}}-\lambda_{2^{\nu-1}}\right)
 E_{2^\nu}(\tilde{f})_p
 \le C
 \sum\limits_{\nu=2^N}^{\infty} \frac{\lambda_{\nu}}{\nu}
\omega_{\alpha+r}\left( f, \frac{1}{\nu}\right)_p,
\end{eqnarray*}
and (\ref{23}) follows.

Repeating the arguments  which were used in (\ref{21}), we estimate
 $n^{-\alpha} \left\|V_n^{(\alpha)}(f^{(\lambda,\beta)})\right\|_p$.
Using Lemma \ref{my2} and inequalities (\ref{11}) and (\ref{20}), we write
\begin{equation}  \label{24}
\left\| V_n^{(\alpha)}(f^{(\lambda,\beta)}) \right\|_p \le
C(\beta, \alpha, r) \left(  n^\alpha \lambda_{n}
\omega_{\alpha+r}\left( f, \frac{1}{n}\right)_p +
\sum\limits_{\nu=1}^{n} \left(\frac{ \lambda_{\nu}}{\nu^r} -
\frac{ \lambda_{\nu+1}}{(\nu+1)^r} \right) \nu^{\alpha+r}
\omega_{\alpha+r}\left( f, \frac{1}{\nu}\right)_p
 \right).
\end{equation}
Collecting (\ref{23}) and (\ref{24}) and using the condition in the right-hand side of (\ref{6}),
 we get $f\in W_p^{ \lambda,\beta} H_\alpha[\varphi].$

Let us prove (\ref{7}). Let $f\in W_p^{ \lambda,\beta}.$ To establish
 $f\in H_{\alpha+\rho}^p(\omega)$, we will first estimate $E_n(f)_p$ for the case $\sin \frac{\pi\beta}{2}\ne 0$.
Then repeating the outline of the proof of estimate (\ref{21}), and taking into account (\ref{20}), we get ( $2^m\le n<2^{m+1}$ )
\begin{eqnarray*}
E_n(f)_p &\le& C \sum\limits_{\nu=m}^{\infty}
\frac{1}{\lambda_{2^\nu}} \left\| \left(V_{[2^{\nu-1}]}-
V_{2^\nu} \right) (f^{(\lambda,\beta)}) \right\|_p \nonumber
\\
&\le&
 C
E_{[2^{m-1}]}(f^{(\lambda,\beta)})_p \sum\limits_{\nu=m}^{\infty}
\frac{1}{\lambda_{2^\nu}} \le C(\rho, p)
   \frac{\|f^{(\lambda,\beta)}\|_p }{\lambda_n}.
\end{eqnarray*}
If $\sin \frac{\pi\beta}{2}= 0$, then it is easy to see that
$$E_n(f)_p \le
\frac{ C}{\lambda_n} E_{n}(f^{(\lambda,\beta)})_p \le
   \frac{C}{\lambda_n}\|f^{(\lambda,\beta)}\|_p . $$
Hence, substituting the obtained bound for $E_n(f)_p$ into inequality (\ref{a3}) and using the fact that
$n^\rho\lambda_n^{-1}\uparrow$ ($n\uparrow$), we obtain
$$
\omega_{\alpha+\rho}\left( f, \frac{1}{n}\right)_p
 \le C(\alpha,\rho) \frac{1}{n^{\alpha+\rho}}
\sum\limits_{\nu=0}^{n} \nu^{\alpha+\rho-1} E_{\nu-1}(f)_p \le
 \frac{C(\alpha,\rho)}{\lambda_n} = O
 \left[\omega\left(\frac{1}{n} \right)\right],$$
i.e., $f\in H_{\alpha+\rho}^p(\omega).$
This completes the proof of the sufficiency part in
(\ref{7}).

Let the right-hand side part of (\ref{8}) hold true and $f\in H^p_{\alpha+\rho}[\omega]$. First let us prove
that
\begin{eqnarray}\label{vallep2} \left\| V_{n}^{(\alpha+r)} (f) \right\|_p
\le C(\alpha,r) \frac{n^r}{\lambda_n} \left\| V_{2n}^{(\alpha)}
(f^{(\lambda,\beta)}) \right\|_p.
\end{eqnarray}
If $\beta=\rho+2m$, and therefore, $r=\rho$, then
$$V_{n}^{(\alpha+r)} (f) = V_{n}^{(\alpha)}
\Bigl(\{\frac{\nu^r}{\lambda_\nu}\}, f^{(\lambda,\beta)}\Bigr)$$
and, by Lemma \ref{my2}
\begin{eqnarray*}
\left\| V_{n}^{(\alpha+r)} (f) \right\|_p \le C(\alpha,r)
\frac{n^r}{\lambda_n}
 \left\| V_{n}^{(\alpha)} (f^{(\lambda,\beta)}) \right\|_p
\le C(\alpha,r) \frac{n^r}{\lambda_n}
  \left\| V_{2n}^{(\alpha)} (f^{(\lambda,\beta)}) \right\|_p.
\end{eqnarray*}
If $\beta\ne\rho+2m$, and therefore, $r>\rho$, then by Lemma \ref{my222}, we write
\begin{eqnarray*}
\left\| V_{2^n}^{(\alpha+r)} (f) \right\|_p &\le&
\sum\limits_{\nu=1}^{n} \left\| V_{2^\nu}^{(\alpha+r)} (f) -
V_{2^{\nu-1}}^{(\alpha+r)} (f)\right\|_p +
 \left\| V_{1}^{(\alpha+r)} (f) \right\|_p
\\
&\le& C \sum\limits_{\nu=1}^{n} \frac{2^{\nu r}}{\lambda_{2^\nu}}
\left\| V_{2^\nu}^{(\alpha)} (f^{(\lambda,\beta)}) -
V_{2^{\nu-1}}^{(\alpha)} (f^{(\lambda,\beta)})\right\|_p +
\frac{1}{\lambda_{1}}
 \left\| V_{1}^{(\alpha)} (f^{(\lambda,\beta)}) \right\|_p
\\
&\le& C \left\| V_{2^{n+1}}^{(\alpha)}
(f^{(\lambda,\beta)})\right\|_p \left(\sum\limits_{\nu=1}^{n}
\frac{2^{\nu r}}{\lambda_{2^\nu}} + \frac{1}{\lambda_{1}}\right)
\\ &\le& C(\alpha,r) \frac{2^{n r}}{\lambda_{2^n}}
  \left\| V_{2^{n+1}}^{(\alpha)} (f^{(\lambda,\beta)}) \right\|_p.
\end{eqnarray*}
 Thus, by (\ref{vallep2}) and the estimate
$$
E_n(f)_p \le \frac{C}{\lambda_n}E_{\left[\frac
n4\right]}(f^{(\lambda,\beta)})_p
$$
we can write
\begin{eqnarray*}
\omega_{\alpha+r}\left( f, \frac{1}{n}\right)_p &\le& C(\alpha,r)
\left( {n}^{-(\alpha+r)} \left\| V_{n}^{(\alpha+r)} (f)
\right\|_p +
E_n(f)_p \right)\\
&\le& \frac{C(\alpha,r)}{\lambda_n} \left( {n}^{-\alpha} \left\|
V_{n}^{(\alpha)} (f^{(\lambda,\beta)}) \right\|_p
  +   E_{\left[\frac n4\right]}(f^{(\lambda,\beta)})_p
\right) \\
&\le& \frac{C(\alpha,r)}{\lambda_n} \omega_{\alpha}\left(
f^{(\lambda,\beta)}, \frac{1}{n}\right)_p = O \left[
 \frac{\varphi\left(\frac{1}{n} \right)}{\lambda_{n}}\right]=
O \left[\omega\left(\frac{1}{n} \right)\right].
\end{eqnarray*}
Thus, the condition
$\frac{\varphi\left(\frac{1}{n} \right)}{\lambda_{n}} = O
\left[\omega\left(\frac{1}{n} \right)\right]$ is sufficient for the  embedding
 $W_p^{\lambda,\beta}H_{\alpha}[\varphi] \subset H_{\alpha+r}^p[\omega]$.

\section{Proof of necessity in Theorem 1.}

We define the trigonometric polynomials $\tau_{n+1} (x)$:
$$
\tau_{n+1} (x) = \sum\limits_{j=1}^{n+1}\alpha_j^{n} \sin j x,
\qquad\mbox{where} \quad \alpha_j^{n}=
\begin{cases} \frac{j}{n+2}, &\text{ $1\le j \le \frac{n+2}{2}$}
\\
1- \frac{j}{n+2}, &\text{ $ \frac{n+2}{2}\le j \le {n+1}$}.
\end{cases}
$$

We will use the following lemmas, as well as Lemmas \ref{my}-\ref{4l22}.
\begin{lemma} \label{lemmagejt} \textnormal{(\cite{tel})}.
Let  series (\ref{fourier}) be the Fourier series of a function $f(x)\in L_1$. Then
$$E_n(f)_1 \ge C \left|
\sum\limits_{\nu=n+1}^\infty \frac{b_\nu}{\nu}\right|.$$
\end{lemma}

\begin{lemma} \textnormal{(\cite[V.1, p. 345; V.2, p.
198]{zygm})}. \label{stlak} Let $1 \le p < \infty$.
\\
\textnormal{(a)} If
the series $\sum\limits_{\nu=1}^{\infty}(a_\nu\cos 2^\nu x+b_\nu\sin
2^\nu x)$ is the Fourier series of a function $f(x)\in L_p$, then
\begin{equation*}
\left\{\sum\limits_{\nu=1}^{\infty} \left( a^2_\nu +
b^2_\nu\right) \right\}^\frac{1}{2} \le C \|f\|_p.
\end{equation*}
\textnormal{(b)} Let $a_n, b_n (n\in \mathbf{N})$ be real
numbers such that
$\sum\limits_{\nu=1}^{\infty} \left( a^2_\nu + b^2_\nu\right)<\infty.$ Then
the series
$\sum\limits_{\nu=1}^{\infty}(a_\nu\cos 2^\nu x +b_\nu\sin 2^\nu x)$
is the Fourier series of a function $f(x)\in L_p$,
and at the same time
\begin{equation*}
 \|f\|_p\le C
\left\{\sum\limits_{\nu=1}^{\infty} \left( a^2_\nu +
b^2_\nu\right) \right\}^\frac{1}{2}.
\end{equation*}
\end{lemma}

\begin{lemma} \textnormal{(\cite[Ch.11,\S 12]{baritext})}.\label{stlak1}
If the series $\sum\limits_{\nu=1}^{\infty}(a_\nu\cos 2^\nu
x+b_\nu\sin 2^\nu x)$,  $a_\nu, b_\nu \ge 0$ is the Fourier series of a function  $f(x)\in L_\infty$, then
$$C_1 \sum\limits_{\xi=n}^{\infty}
(a_\xi+b_\xi) \le  E_{2^n-1}(f)_\infty \le C_2
\sum\limits_{\xi=n}^{\infty} (a_\xi+b_\xi) .$$
\end{lemma}

We will use the following definitions. Let $\omega(\cdot)\in\Phi_{\alpha}$.
\\
A sequence $\psi$ is called  {\em $Q_{\alpha,\, \theta}(\omega)$-sequence} if
\begin{eqnarray}
0&<& \psi_n \,\,\le\,\, n^{\alpha} \omega \left(
\frac{1}{n}\right),\quad \psi_n \uparrow \,(n\uparrow)\label{zv2}
\\
C_1 \,\omega \left( \frac{1}{n}\right) &\le& \left\{
\sum\limits_{\nu=n}^\infty \nu^{-\alpha\theta-1}  \psi_\nu^\theta
\right\}^\frac{1}{\theta} \le\, C_2 \,\omega \left(
\frac{1}{n}\right). \label{27}
\end{eqnarray}
A sequence $\varepsilon$ is called  {\em $q_{\alpha,\, \theta}(\omega)$-sequence} if
\begin{eqnarray}
0&<& \varepsilon_n \,\,\le\,\, \omega \left(
\frac{1}{n+1}\right),\quad \varepsilon_n \downarrow
\,(n\uparrow)\label{32}
\\
C_1 \,\omega \left( \frac{1}{n+1}\right) &\le&
 \left\{ (n+1)^{-\alpha\theta}
\sum\limits_{\nu=1}^{n+1} \nu^{\alpha\theta-1}
\varepsilon_\nu^\theta \right\}^\frac{1}{\theta} \le\,
 C_2 \,\omega \left( \frac{1}{n+1}\right). \label{33}
\end{eqnarray}

\begin{center}{\bf Necessity in (\ref{1}) - (\ref{8}).}
\end{center}
We prove the necessity part by constructing corresponding examples. The proof consists of eight steps.

\textbf{I.} $1<p<\infty$. {\sf {Step 1.} } Let us show the necessity part in (\ref{1}).
\\
Let $\omega(\cdot)\in\Phi_{\alpha+\rho}$ and $\theta=\min(2, p)$.
We will construct a $Q_{\alpha+\rho,\, \theta}(\omega)$-sequence $\psi$.
\\
Assume that integers $1=n_1<n_2<\cdots<n_s$ are chosen.
Then as $n_{s+1}$ we take the minimum number $N>n_s$ such that
$$
 \omega \left( \frac{1}{N}\right)<
\frac{1}{2} \,\,\omega \left( \frac{1}{n_s}\right) \le  \omega
\left( \frac{1}{N-1}\right).
$$
  We set
$$
\psi_n= \begin{cases} n_s^{\rho+\alpha} \omega \left(
\frac{1}{n_s}\right), &\text{if} \quad n_s\le n <n_{s+1},
\quad
s=1,2,\cdots;
\\
0,
&\text{if}
\quad n=0.
\end{cases}
$$
It is easy to see that this sequence is what we need.

Let $H^p_{\alpha+\rho}[\omega] \subset W_p^{ \lambda,\beta}$  and
let the series in (\ref{1}) be divergent.
By means of properties of sequence $\left\{\psi_n\right\}$, we have
\begin{eqnarray*}
\infty= \sum\limits_{\nu=1}^\infty \left( \lambda_{\nu+1}^\theta
- \lambda_{\nu}^\theta \right) \omega^\theta \left(
\frac{1}{\nu}\right) &\le&
 C(\alpha,\rho,\theta)
\sum\limits_{\nu=1}^\infty \left( \lambda_{\nu+1}^\theta -
\lambda_{\nu}^\theta \right) \sum\limits_{m=\nu}^\infty
m^{-(\alpha+\rho)\theta-1} \psi_m^\theta
\\
 &\le&
 C(\alpha,\rho,\theta)
\sum\limits_{\nu=1}^\infty \lambda_{\nu}^\theta
\nu^{-(\alpha+\rho)\theta-1} \psi_\nu^\theta.
\end{eqnarray*}
\\
{\sf Step 1(a): $2\le p <\infty$}.
    We consider the series
\begin{equation}\label{28}
\sum\limits_{\nu=1}^\infty 2^{-\nu(\alpha+\rho)} \left(
\psi_{2^\nu}^2- \psi_{2^{\nu-1}}^2 \right)^\frac{1}{2} \cos 2^\nu
x.
\end{equation}
Since
\begin{eqnarray}
\label{zv3} \sum\limits_{\nu=1}^\infty 2^{-2\nu(\alpha+\rho)}
\left( \psi_{2^\nu}^2- \psi_{2^{\nu-1}}^2 \right) &\le&
\sum\limits_{\nu=1}^\infty
 \left(
\psi_{2^\nu}^2- \psi_{2^{\nu-1}}^2 \right)
\sum\limits_{\xi=\nu}^\infty 2^{-2\xi(\alpha+\rho)}
\nonumber\\
&\le& \sum\limits_{\nu=1}^\infty 2^{-2\nu(\alpha+\rho)}
\psi_{2^\nu}^2 \le  C \omega^2(1),
\end{eqnarray}
then, by Zygmund's Lemma \ref{stlak}, series (\ref{28}) is the
Fourier series of a function $f_1(x)\in L_p$. Applying Lemmas \ref{my}
and \ref{stlak}, we get
$$
C(\alpha,\rho)
 \omega_{\alpha+\rho}\Bigl(f_1,\frac {1}{2^{n}}\Bigr)_p \le  2^{-n(\alpha+\rho)}
\left(
 \sum\limits_{\nu=1}^{n} a_\nu^2 2^{2(\alpha+\rho) \nu}\right)^\frac{1}{2}+
\left(\sum\limits_{\nu={n+1}}^{\infty} a_\nu^2
\right)^{\frac{1}{2}}=: I_1+I_2,$$ where
$a_\nu=2^{-(\alpha+\rho)}\left( \psi_{2^\nu}^2-
\psi_{2^{\nu-1}}^2 \right)^\frac{1}{2}.$
Similarly to estimate (\ref{zv3}), by means of (\ref{zv2}), we get
$$
I_1\le
2^{-n(\alpha+\rho)} \psi_{2^n} \le
 \omega\Bigl(\frac{1}{2^n}\Bigr),
$$
and by means of (\ref{27}),
$$
I_2\le C(\alpha,\rho)
\left( \sum\limits_{\nu=n+1}^\infty  2^{-2\nu(\alpha+\rho)}
\psi_{2^\nu}^2 \right)^\frac{1}{2} \le C(\alpha,\rho)
\omega\Bigl(\frac{1}{2^n}\Bigr).
$$

Thus, $f_1(x)\in H^p_{\alpha+\rho}[\omega].$ Then from our
assumption, $f_1(x)\in W_p^{\lambda,\beta}$.
On the other hand,
$$\left\|f_1^{(\lambda,\beta)}\right\|_p \ge
 C(\alpha,\rho,\theta)
\left( \sum\limits_{\nu=1}^\infty \lambda_{\nu}^2
\nu^{-2(\alpha+\rho)-1} \psi_\nu^2 \right)^\frac{1}{2}\, =
\,\infty.$$
This contradiction proves the convergence of series
in (\ref{1}).
\\
{\sf Step 1(b): $1<p\le 2$}.
 Consider series\footnote{Series of this type was considered in \cite{mtiman}.}
\begin{equation}\label{29}
\psi_{1} \cos  x + \sum\limits_{\nu=1}^\infty
2^{-\nu(\alpha+\rho)} 2^{\nu(\frac{1}{p}-1)} \left(
\psi_{2^\nu}^p- \psi_{2^{\nu-1}}^p \right)^\frac{1}{p}
\sum\limits_{\mu=2^{\nu-1}+1}^{2^{\nu}} \cos \mu x.
\end{equation}
Using the Jensen inequality\,
 $\left( \sum\limits_{n=1}^{\infty} a_n^\alpha \right)^{1/\alpha} \le
\left( \sum\limits_{n=1}^{\infty} a_n^\beta \right)^{1/\beta}$
($a_n\ge 0$ and $\,0 < \beta \le \alpha<\infty$), we write
\begin{eqnarray*}
& & \int\limits_0^{2\pi} \left[ \sum\limits_{\nu=1}^\infty \left(
2^{-\nu(\alpha+\rho)} 2^{\nu(\frac{1}{p}-1)} \left(
\psi_{2^\nu}^p- \psi_{2^{\nu-1}}^p \right)^\frac{1}{p}
\sum\limits_{\mu=2^{\nu-1}+1}^{2^{\nu}} \cos \mu x \right)^2
\right]^\frac{p}{2} dx
\\
&\le& \int\limits_0^{2\pi} \left[ \sum\limits_{\nu=1}^\infty
2^{-\nu p(\alpha+\rho)} 2^{\nu(1-p)} \left( \psi_{2^\nu}^p-
\psi_{2^{\nu-1}}^p \right) \left|
\sum\limits_{\mu=2^{\nu-1}+1}^{2^{\nu}} \cos \mu x \right|^p
\,\right] dx
\\
&\le& C(p) \sum\limits_{\nu=1}^\infty \left( \psi_{2^\nu}^p-
\psi_{2^{\nu-1}}^p \right) 2^{-\nu p(\alpha+\rho)} \le C(p)
\omega^p(1),
\end{eqnarray*}
because of $\,\,C_1(p) 2^{\nu(p-1)}\le  \left\|
\sum\limits_{\mu=2^{\nu-1}+1}^{2^{\nu}} \cos \mu x \right\|^p_p
\le C_2(p)2^{\nu(p-1)}$.

By the Littlewood-Paley theorem (see \cite[Vol. 2, p. 349]{zygm}), there exists a function $f_2\in L_p$ with the Fourier series (\ref{29}).
One can easily check that $f_2 \in H^p_{\alpha+\rho}[\omega].$
By our assumption, $f_2(x)\in W_p^{\lambda,\beta}$.
On the other hand, Paley's theorem on Fourier coefficients \cite[V.2, p.
182]{zygm} implies that for $f_2\in L_p$
\begin{eqnarray*}
\left\|f_2^{(\lambda,\beta)}\right\|_p^p &\ge&
 C(p)
\sum\limits_{\nu=1}^\infty 2^{-\nu p(\alpha+\rho)} 2^{\nu(1-p)}
\left( \psi_{2^\nu}^p- \psi_{2^{\nu-1}}^p \right)
\sum\limits_{\mu=2^{\nu-1}+1}^{2^{\nu}} \lambda_\mu^p \mu^{p-2}
\\
&\ge&
 C(\alpha,\rho, p)
\sum\limits_{\nu=1}^\infty \left( \psi_{2^\nu}^p-
\psi_{2^{\nu-1}}^p \right) \sum\limits_{\xi=\nu}^\infty \left(
2^{-\xi(\alpha+\rho)p}\lambda_{2^\xi}^p -
2^{-(\xi+1)(\alpha+\rho)p}\lambda_{2^{\xi+1}}^p \right)
\\
&\ge&
 C_1(\alpha,\rho,p)
\sum\limits_{\nu=1}^\infty \psi_\nu^p \lambda_{\nu}^p
\nu^{-p(\alpha+\rho)-1} -
 C_2(\alpha,\rho,p)
 \psi_1^p \lambda_{2}^p
2^{-p(\alpha+\rho)-1}\, = \,\infty.
\end{eqnarray*}
This contradiction shows that the series in the right-hand side of (\ref{1}) converges. This completes
the proof of the necessity part of (\ref{1}).

{\sf {Step 2.} } Let us prove the necessity in (\ref{2}) for the case $2\le p <\infty.$
\\
We notice that using Lemmas \ref{my} and \ref{stlak}, we get for
$$
f(x)\sim \sum\limits_{\nu=1}^{\infty} (a_\nu \cos 2^\nu x +
b_\nu \sin 2^\nu x)
$$
the following relation
\begin{equation}\label{mod-cF}
\omega_\alpha\left(f, \frac{1}{2^m}\right)_p \asymp \left(
2^{-2m\alpha} \sum\limits_{\nu=1}^{m} (a_\nu^2 + b_\nu^2)
2^{2\nu\alpha}\right)^\frac{1}{2}+
\left(\sum\limits_{\nu=m+1}^{\infty} (a_\nu^2 +
b_\nu^2)\right)^\frac{1}{2}.
\end{equation}
Let $\omega(\cdot)\in\Phi_{\alpha+\rho}$. Then one can construct\footnote{See, for example, \cite{temirgaliev}, \cite{gejt2}.} a sequence $\varepsilon$ such that it is a $q_{\alpha+\rho, \,\theta}(\omega)$-sequence.
In this case, we consider
\begin{equation}\label{34}
\varepsilon_0+ \left( \varepsilon_{1}^2 -
\varepsilon_{2}^2\right)^\frac{1}{2} \cos  x +
\sum\limits_{\nu=1}^\infty \left( \varepsilon_{2^{\nu}}^2 -
\varepsilon_{2^{\nu+1}}^2\right)^\frac{1}{2} \cos 2^{\nu} x.
\end{equation}
Repeating the argument used for series (\ref{28}), we obtain that series (\ref{34}) is the Fourier
series of a function $f_3\in L_p$. Since $E_{2^{n-1}}(f_3)_p\le C(p)\varepsilon_{2^{n}}$, then by
(\ref{a3}) and (\ref{33}) we have $f_3\in H_{\alpha+\rho}^p[\omega]$. We define
$f_{13}:=f_{1}+f_{3}.$
Then  $f_{13}\in H_{\alpha+\rho}^p[\omega]\subset W_p^{ \lambda,\beta} H_\alpha[\varphi]$.
 It is easy to see from (\ref{mod-cF}) that
 \begin{equation}\label{99}
C(\alpha, \beta) \,
\omega_{\alpha}\Bigl(f_{13}^{(\lambda,\beta)},
\frac{1}{n+1}\Bigr)_p \ge
\omega_{\alpha}\Bigl(f_{1}^{(\lambda,\beta)},
\frac{1}{n+1}\Bigr)_p+
\omega_{\alpha}\Bigl(f_{3}^{(\lambda,\beta)},
\frac{1}{n+1}\Bigr)_p.
 \end{equation}

Let us estimate $\omega_{\alpha}\Bigl(f_{1}^{(\lambda,\beta)},
\frac{1}{n+1}\Bigr)_p$.
 Applying (\ref{mod-cF}) and using the properties of the sequence  $\{\psi_{\nu}\}$,
 we write ($2^m\le n+1< 2^{m+1}$)
\begin{eqnarray}
\omega^2_{\alpha}\left(f_1^{(\lambda,\beta)},
\frac{1}{n+1}\right)_p &\ge&  C(\alpha, \beta)
\sum\limits_{\nu=m}^\infty \psi_{2^\nu}^2 2^{-2\nu(r+\alpha)}
\left(
 \sum\limits_{k=m}^{\nu} (\lambda_{2^k}^2 - \lambda_{2^{k-1}}^2) +
 \lambda_{2^{m-1}}^2\right)
 \nonumber
\\
&\ge& C(\alpha, \beta) \left( \lambda_{2^m}^2
\omega^2\left(\frac{1}{2^m}\right) + \sum\limits_{k=m}^\infty
(\lambda_{2^k}^2 - \lambda_{2^{k-1}}^2)
\omega^2\left(\frac{1}{2^k}\right)\right)
 \nonumber
 \\
&\ge& C(\alpha, \beta) \left(
 \lambda_{n+1}^2 \omega^2\left(\frac{1}{n+1}\right) + \sum\limits_{\nu={n+2}}^\infty
\left(\lambda_{\nu+1}^2 -\lambda_{\nu}^2\right)
\omega^2\left(\frac{1}{\nu}\right) \right) .
 \label{f1}
\end{eqnarray}

Let us proceed to the estimation of $\omega_{\alpha}\Bigl(f_{3}^{(\lambda,\beta)},
\frac{1}{n+1}\Bigr)_p.$ Using (\ref{mod-cF}), we have
($2^m\le n+1< 2^{m+1}$):
\begin{eqnarray} \label{35}
\omega_\alpha^2\left(f_3^{(\lambda,\beta)},
\frac{1}{n+1}\right)_p &\ge& C(\alpha, \beta) 2^{-2 m \alpha}
 \sum\limits_{\nu=0}^m 2^{2 \nu \alpha} \lambda_{2 \nu}^2
\left( \varepsilon_{2^{\nu}}^2 - \varepsilon_{2^{\nu+1}}^2\right)
\nonumber
\\
&\ge& C_1(\alpha, \beta) 2^{-2 m \alpha}
 \sum\limits_{\nu=0}^m 2^{2 \nu \alpha} \lambda_{2^\nu}^2
\varepsilon_{2^{\nu}}^2 - C_2(\alpha,p)    \lambda_{2^{m+1}}^2
\varepsilon_{2^{m+1}}^2.
\end{eqnarray}
Then Jackson's inequality implies
\begin{eqnarray} \label{36}
\omega_\alpha^2\left(f_3^{(\lambda,\beta)},
\frac{1}{n+1}\right)_p &\ge& C(\alpha)
E_{2^m-1}^2\left(f_3^{(\lambda,\beta)} \right)_p \nonumber \\
&\ge& C(\alpha,
\beta)
 \sum\limits_{\nu=m}^\infty  \lambda_{2 \nu}^2
\left( \varepsilon_{2^{\nu}}^2 - \varepsilon_{2^{\nu+1}}^2\right) \nonumber \\
&\ge& C(\alpha, \beta)
 \lambda_{2^m}^2 \varepsilon_{2^{m}}^2.
\end{eqnarray}
Applying estimates (\ref{35}) and (\ref{36}), we get
\begin{equation} \label{37}
\omega_\alpha\left(f_3^{(\lambda,\beta)}, \frac{1}{n+1}\right)_p
\ge C(\alpha, \beta) \left((n+1)^{-2\alpha}
 \sum\limits_{\nu=1}^{n+1}
  \lambda_{\nu}^2
\nu^{2\alpha-1} \varepsilon^2_{\nu}\right)^\frac{1}{2}.
\end{equation}
Further, (\ref{33}) and $\nu^{-\rho}\lambda_\nu \downarrow$ allow us to write the estimates
\begin{eqnarray}
  \lambda_{n+1}^2 \omega^2 \left( \frac{1}{n+1}\right)
&+& (n+1)^{-2\alpha} \sum\limits_{\nu=1}^{n+1}
\nu^{2(\rho+\alpha)} \omega^2 \left( \frac{1}{\nu}\right) \left(
  \lambda_{\nu}^2 \nu^{-2\rho}-
  \lambda_{\nu+1}^2 (\nu+1)^{-2\rho}
\right)\nonumber
\\
&\le& C(\alpha, \rho) (n+1)^{-2\alpha} \sum\limits_{\nu=1}^{n+1}
\nu^{2\alpha-1}
 \lambda_{\nu}^2 \varepsilon_{\nu}^2.
\label{f3}
\end{eqnarray}
Combining estimates   (\ref{99}), (\ref{f1}), (\ref{37}), (\ref{f3}), and
 $\omega_{\alpha}\Bigl(f_{13}^{(\lambda,\beta)},
\frac{1}{n+1}\Bigr)_p =
O\left[\varphi\left(\frac{1}{n}\right)\right]$,
we obtain the condition in the right-hand side part of
(\ref{2}).

{\sf {Step 3.} } Here we show the necessity in (\ref{2}) for the case of $1<p<2.$
In this case the proof  is similar to the proof of the case $2\le p <\infty$.
The only difference is that we use Paley's theorem on Fourier coefficients instead of Zygmund's Lemma \ref{stlak}.
In this case we consider the sum of $f_{2}(x)$ and the
following function
\begin{equation}\label{38}
\varepsilon_0+ \left( \varepsilon_{1}^p -
\varepsilon_{2}^p\right)^\frac{1}{p} \cos  x +
\sum\limits_{\nu=0}^\infty 2^{\nu\left(\frac{1}{p}-1\right)}
\left( \varepsilon_{2^{\nu+1}}^p -
\varepsilon_{2^{\nu+2}}^p\right)^\frac{1}{p}
\sum\limits_{\mu=2^{\nu}+1}^{2^{\nu+1}} \cos \mu x.
\end{equation}

{\sf {Step 4.}} To prove the necessity in (\ref{3}) and (\ref{7}), we consider the general case
of $1\le p\le \infty.$
Let $\Phi$ be the class of all decreasing null-sequences. It is clear that
\begin{equation*}\label{conditions}
\frac{1}{\lambda_n}= O \left[ \omega\left(\frac{1}{n}\right)
\right] \qquad \Longleftrightarrow \qquad
\forall\,\,\gamma=\{\gamma_n\}\in \Phi\qquad
\frac{\gamma_n}{\lambda_n}= O \left[
\omega\left(\frac{1}{n}\right) \right].
\end{equation*}
Let us assume that $\frac{\gamma_n}{\lambda_n}= O
\left[\omega\left(\frac{1}{n}\right)\right]$ does not hold for all
$\,\gamma\in \Phi$ and $W^{\lambda,\beta}_p \subset
H^p_{\alpha+\rho}[\omega]$.  Then there exist
$\gamma=\{\gamma_n\}\in \Phi$ and $\{C_{n}\uparrow \infty\}$ such that $\frac{\gamma_{m_n}}{\lambda_{m_n}}\ge C_{n}
\omega\left(\frac{1}{m_n}\right).$
Further, we choose a subsequence  $\left\{m_{n_k} \right\}$ such that
$\frac{m_{n_{k+1}}}{m_{n_k}}\ge 2$ and $\gamma_{m_{n_k}}\le 2^{-k}.$
Consider the series
\begin{equation}\label{19vs}
\sum\limits_{k=0}^\infty
\frac{\gamma_{m_{n_k}}}{\lambda_{m_{n_k}}} \cos ( {m_{n_k}}+1)x.
\end{equation}
Since $\sum\limits_{k=0}^\infty
\frac{\gamma_{m_{n_k}}}{\lambda_{m_{n_k}}} \le
\frac{1}{\lambda_{m_{n_0}}} \sum\limits_{k=0}^\infty
\frac{1}{2^k}<\infty,$
 there exists a function $f_4\in L_p$ with the Fourier series (\ref{19vs}).
 Because of
 $\sum\limits_{k=0}^\infty {\gamma_{m_{n_k}}} \le
\sum\limits_{k=0}^\infty \frac{1}{2^k}<\infty,$ we have
$f_4^{(\lambda,\beta)}\in L_p$, i.e., $f_4\in W_p^{\lambda,\beta}$.

On the other hand, using (\ref{11}) and
 $E_{n-1}(f)_p \ge C \left(|a_n|+|b_n|\right)$,
\begin{eqnarray*}
\omega_{\alpha+\rho} \left(f_4, \frac{1}{m_{n_k}} \right)_p\ge
C(\alpha, \rho) E_{m_{n_k}}(f_{4})_p \ge C(\alpha, \rho)
\frac{\gamma_{m_{n_k}}}{\lambda_{m_{n_k}}} \ge C(\alpha, \rho)
C_{n_k} \omega\left(\frac{1}{m_{n_k}}\right),
\end{eqnarray*}
i.e.,  $f_4\notin H^p_{\alpha+\rho}[\omega].$
This contradiction proves that the condition
 $\frac{1}{\lambda_n}= O
\left[\omega\left(\frac{1}{n}\right)\right]$ is necessary for
$W_p^{\lambda,\beta} \subset H^p_{\alpha+\rho}[\omega]$.

{\sf {Step 5.}} To prove the necessity in (\ref{4}) and (\ref{8}), we verify that for any
$\rho>0$ and $1\le p\le \infty$,
\begin{equation}\label{emb1}
W_p^{ \lambda,\beta} H_\alpha[\varphi]
 \subset H_{\alpha+\rho}^p[\omega]
\quad \Longrightarrow \quad \frac{\varphi\left(\frac{1}{n}
\right)}{\lambda_{n}}= O \left[\omega\left(\frac{1}{n}
\right)\right].
\end{equation}
First, we remark that
\begin{equation*}
\frac{\varphi\left(\frac{1}{n} \right)}{\lambda_{n}}= O
\left[\omega\left(\frac{1}{n} \right)\right] \qquad
\Longleftrightarrow \qquad \forall\,\,\gamma=\{\gamma_n\}\in
\Phi\qquad \frac{\gamma_n \varphi\left(\frac{1}{n}
\right)}{\lambda_{n}}= O \left[\omega\left(\frac{1}{n}
\right)\right].
\end{equation*}
We assume that the relation in the right-hand side of (\ref{emb1}) does not hold.
Then there exist $\gamma=\{\gamma_n\}\in \Phi$ and $\{C_{n}\uparrow \infty\}$ such that
$\frac{\gamma_{m_n}\varphi\left(\frac{1}{m_n}\right)}{\lambda_{m_n}}\ge
C_{n} \omega\left(\frac{1}{m_n}\right).$
We choose a subsequence   $\left\{m_{n_k} \right\}$ such that
$\frac{m_{n_{k+1}}}{m_{n_k}}\ge 2$ ш $\gamma_{m_{n_k}}\le 2^{-k}.$
Because of $\sum\limits_{k=0}^\infty
\gamma_{m_{n_k}}\varphi\left(\frac{1}{m_n}\right)
 \le
\varphi\left(\frac{1}{m_0}\right) \sum\limits_{k=0}^\infty
\frac{1}{2^k}<\infty,$ there exists a
function $f_5\in L_p$ with the Fourier series
\begin{equation}\label{19vsvs}
\sum\limits_{k=0}^\infty
{\gamma_{m_{n_k}}}\varphi\left(\frac{1}{m_{n_k}}\right) \cos (
{m_{n_k}}+1)x.
\end{equation}
For  $m_{n_k}\le n < m_{n_{k+1}}$ using  Lemmas \ref{my} and \ref{stlak1}, we have
\begin{eqnarray*}
\omega_{\alpha} \left(f_5, \frac{1}{n} \right)_p
 &\le&
 C \omega_{\alpha} \left(f_5, \frac{1}{n} \right)_\infty
 \le C \left( n^{-\alpha}
 \sum \limits_{s=0}^k
{\gamma_{m_{n_s}}}\varphi\left(\frac{1}{m_{n_s}}\right)
{m^{\alpha}_{n_s}} + \sum
\limits_{s=k+1}^\infty {\gamma_{m_{n_s}}}\varphi\left(\frac{1}{m_{n_s}}\right) \right) \\
 &\le&
C\left( \varphi\left(\frac{1}{n}\right)
 \sum \limits_{s=0}^k
{\gamma_{m_{n_s}}} + \varphi\left(\frac{1}{n}\right)
 \sum \limits_{s=k+1}^\infty
{\gamma_{m_{n_s}}} \right) \le  C \varphi\left(\frac{1}{n}\right).
\end{eqnarray*}
Then $f_5\in H^p_{\alpha}[\varphi]$, i.e., setting
$\frac{1}{\lambda} := \left\{ \frac{1}{\lambda_1},
\frac{1}{\lambda_2}, \frac{1}{\lambda_3}, \cdots \right\}$, we have
$f_5^{(\frac{1}{\lambda}, - \beta)} \in W_p^{ \lambda,\beta}
H_\alpha[\varphi].$

On the other hand,
\begin{eqnarray*}
\omega_{\alpha+\rho} \left(f_5^{(\frac{1}{\lambda}, - \beta)},
\frac{1}{m_{n_k}} \right)_p\ge C
E_{m_{n_k}}(f_5^{(\frac{1}{\lambda}, - \beta)})_p \ge C
\frac{\gamma_{m_{n_k}} \varphi\left(\frac{1}{m_{n_k}}\right)}
{\lambda_{m_{n_k}}} \ge C C_{n_k}
\omega\left(\frac{1}{m_{n_k}}\right),
\end{eqnarray*}
i.e., $f_5^{(\frac{1}{\lambda}, - \beta)}\notin
H^p_{\alpha+\rho}[\omega].$ This contradicts our assumption.
The proof of the necessity part in (\ref{3})-(\ref{4}) and (\ref{7})-(\ref{8}) is now complete.

\textbf{II.} $p=1$ or $p=\infty$. {\sf {Step 6.} }
Let us prove the necessity in  (\ref{5}). Let  $H_{\alpha+\rho}^p[\omega] \subset
W_p^{\lambda,\beta}$ and the series in (\ref{5}) be divergent.

{\sf Step 6(a): $\sin \frac{\beta\pi}{2} \ne 0$}.
In this case a divergence of the series in (\ref{5})
is equivalent to the divergence of the series
$\sum\limits_{n=1}^\infty \frac{\lambda_{n}}{n} \omega\left(\frac{1}{n}\right).$

Let $p=1$. We take a $q_{\alpha+\rho, \,1}(\omega)$-sequence $\varepsilon$
and consider the series
\begin{equation}\label{39} \sum\limits_{\nu=1}^\infty
\left(\varepsilon_\nu - \varepsilon_{\nu+1}\right) K_\nu(x).
\end{equation}
This series converges in $L_1$ (see \cite{gejt}) to a function $f_6(x)$ and $E_n(f_6)_1=O\left(\varepsilon_n\right)$.
Applying (\ref{a3}) and (\ref{33}), we get $f_6 \in H_{\alpha+\rho}^p[\omega]
\subset W_p^{\lambda, \beta}$.
One can also rewrite (\ref{39}) in the following form
$$
\sum\limits_{\nu=1}^\infty  a_\nu \cos \nu x, \qquad \mbox{where}
\qquad
 a_\nu=\varepsilon_\nu -\nu \sum\limits_{j=\nu}^\infty
\frac{\varepsilon_j - \varepsilon_{j+1}}{j+1}.
$$
By Lemma \ref{lemmagejt},
\begin{eqnarray*}
\left\| f_{6}^{(\lambda,\beta)} \right\|_1 &\ge& C(\beta)
\sum\limits_{\nu=1}^\infty \frac{\lambda_\nu}{\nu} a_{\nu} =
 C(\beta)
 \left(
\sum\limits_{\nu=1}^\infty
\frac{\lambda_\nu}{\nu}\varepsilon_{\nu} -
\sum\limits_{\nu=1}^\infty \lambda_\nu \sum\limits_{j=\nu}^\infty
\frac{\varepsilon_{j} - \varepsilon_{j+1}}{j+1} \right)
\\
&=& C(\beta) \left( \sum\limits_{\nu=1}^\infty
\frac{\lambda_\nu}{\nu}\varepsilon_{\nu} -
\sum\limits_{\nu=1}^\infty (a_{\nu}-a_{\nu+1}) \lambda_{\nu}
 \right)
\\
&\ge& C_1(\beta)
 \sum\limits_{\nu=1}^\infty
\frac{\lambda_\nu}{\nu}\varepsilon_{\nu} - C_2(\beta) \left(
 \lambda_{1} a_{1} +
 \sum\limits_{n=1}^\infty
 (\lambda_{n+1}-\lambda_{n})a_{n} \right).
 \end{eqnarray*}
Further, using (\ref{32}) and (\ref{33}), we get
\begin{eqnarray*}
 \sum\limits_{\nu=1}^\infty
\frac{\lambda_\nu}{\nu}\varepsilon_{\nu} &\le& C(\rho)
 \sum\limits_{\nu=1}^\infty
\frac{\lambda_\nu}{\nu} \omega\left(\frac{1}{\nu}\right)  \le
C(\rho)
 \sum\limits_{\nu=1}^\infty
{\lambda_\nu} \nu^{-(\alpha+\rho)-1}
 \sum\limits_{m=1}^\nu
m^{\alpha+\rho-1} \varepsilon_{m}
\\
&=& C(\rho)
 \sum\limits_{m=1}^\infty m^{\alpha+\rho-1} \varepsilon_{m}
 \sum\limits_{\nu=m}^\infty
\frac{\lambda_\nu}{\nu^\rho} \nu^{-\alpha-1} \le C(\rho, \alpha)
 \sum\limits_{m=1}^\infty
\frac{\lambda_m}{m}\varepsilon_{m}.
 \end{eqnarray*}
Then
\begin{equation}\label{41}
\left\| f_6^{(\lambda,\beta)} \right\|_1 \ge C_1(\alpha, \rho,
\beta)
 \sum\limits_{\nu=1}^\infty \lambda_\nu \nu^{-1}
\omega\left(\frac{1}{\nu}\right)- C_2(\alpha, \rho, \beta) \left(
\lambda_1 a_1 + \sum\limits_{\nu=1}^\infty (\lambda_{\nu+1} -
\lambda_{\nu}) a_\nu \right).
\end{equation}
On the other hand, using monotonicity of $\{a_\nu\}$ and Lemma \ref{lemmagejt}, we have
\begin{eqnarray}\label{42}
C(\beta, \rho)\left\| f_6^{(\lambda,\beta)} \right\|_1 &\ge&
C(\rho) \left( \lambda_1 a_1 + \sum\limits_{\nu=1}^\infty
\frac{\lambda_{\nu} a_\nu}{\nu} \right) \ge C(\rho) \left(
\lambda_1 a_1 + \sum\limits_{\nu=0}^\infty \lambda_{2^\nu}
a_{2^\nu} \right)
\nonumber\\
&\ge& \lambda_1 a_1 + \sum\limits_{\nu=0}^\infty  a_{2^\nu} (
\lambda_{2^{\nu+1}} - \lambda_{2^{\nu}}) \ge \lambda_1 a_1 +
\sum\limits_{\nu=1}^\infty (\lambda_{\nu+1} - \lambda_{\nu})
a_\nu.
\end{eqnarray}
Using  (\ref{41}) and (\ref{42}), we get
$$
\left\| f_6^{(\lambda,\beta)} \right\|_1 \ge C(\alpha, \rho,
\beta) \sum\limits_{\nu=1}^\infty \lambda_\nu \nu^{-1}
\omega\left(\frac{1}{\nu}\right) = \infty.
$$
The obtained contradiction implies the convergence of series in (\ref{5}).

Let now $p=\infty$. Define the function  (see \cite{bari-sopr})
$f_7(x) = \sum\limits_{\nu=1}^\infty \varepsilon_\nu \nu^{-1} \sin
\nu x,$ where
 $\varepsilon$ is a $q_{\alpha+\rho, \,1}(\omega)$-sequence.
 Then  $E_n(f_7)_\infty \le C\varepsilon_{n+1}$.
 Using (\ref{a3}) and (\ref{33}), we get
$f_7\in H_{\alpha+\rho}^p[\omega]\subset W_p^{\lambda,\beta}$. On the other hand,
\begin{eqnarray*}
\infty = \left\| f_7^{(\lambda,\beta)} \right\|_\infty &\ge&
C(\beta)
 \sum\limits_{\nu=0}^\infty 2^{\nu(\alpha+\rho)}
\varepsilon_{2^\nu} \frac{\lambda_{2^\nu}}{2^{\nu(\alpha+\rho)}}
\\
&\ge& C(\alpha, \rho, \beta) \sum\limits_{\nu=0}^\infty
2^{\nu(\alpha+\rho)}\varepsilon_{2^\nu}
\sum\limits_{m=\nu}^\infty
\frac{\lambda_{2^m}}{2^{m(\alpha+\rho)}}
\\
&\ge&  C(\alpha, \rho, \beta)
 \sum\limits_{\nu=1}^\infty
\lambda_\nu \nu^{-1} \omega\left(\frac{1}{\nu}\right).
\end{eqnarray*}
The obtained contradiction proves the convergence of the series in (\ref{5}).

{\sf Step 6(b): $\sin \frac{\beta\pi}{2} = 0$}.
Let the series in (\ref{5}) be divergent. We consider only the non-trivial case of $\rho>0.$ Let
 $\varepsilon$  be a $q_{\alpha+\rho,\,1}(\omega)$-sequence. By means of the properties
  $\{\lambda_n\}$, we have
\begin{eqnarray}
\sum\limits_{\nu=2}^{\infty} \left(
  \lambda_{\nu+1} - \lambda_{\nu}
\right) \omega \left( \frac{1}{\nu}\right) &\le& C(\alpha, \rho)
\sum\limits_{s=0}^{\infty} \left(\varepsilon_{2^s} -
\varepsilon_{2^{s+1}} \right) \left[ \sum\limits_{m=0}^{s} \right.
2^{m(\alpha+\rho)} \sum\limits_{\nu=m}^{s} 2^{-\nu(\alpha+\rho)}
\left(
  \lambda_{2^\nu} - \lambda_{2^{\nu-1}}  \right)
\nonumber
\\
&+& \left.\sum\limits_{m=0}^{s} 2^{m(\alpha+\rho)}
\sum\limits_{\nu=s+1}^{\infty} 2^{-\nu(\alpha+\rho)} \left(
  \lambda_{2^\nu} - \lambda_{2^{\nu-1}}  \right)
 \right]
\nonumber
\\
&\le& C(\alpha, \rho) \sum\limits_{\nu=0}^{\infty}
\left(\varepsilon_{\nu} - \varepsilon_{\nu+1} \right)
  \lambda_{\nu}.
\label{43}
\end{eqnarray}
Let $p=1$.  Then the series
\begin{equation}\label{44}
\sum\limits_{\nu=1}^\infty \left(\varepsilon_{\nu} -
\varepsilon_{\nu+1} \right) \tau_\nu(x)
\end{equation}
converges (\cite{gejt}) to $f_8\in L_1$ and $E_n(f_8)_1\le C
\varepsilon_{n+1}$. Hence, we have $f_8\in
H_{\alpha+\rho}^p[\omega]\subset W_p^{\lambda,\beta}$.

We write series (\ref{44}) in the following way
$$
\sum\limits_{\nu=1}^\infty  b_\nu \sin \nu x, \qquad \mbox{where}
\qquad
 b_\nu=
\sum\limits_{j=\nu}^{2\nu-2}
 \left( 1 - \frac{\nu}{j+1} \right)
\left(\varepsilon_{\nu} - \varepsilon_{\nu+1} \right) +
\sum\limits_{j=2\nu-1}^\infty
 \frac{\nu}{j+1}
\left(\varepsilon_{\nu} - \varepsilon_{\nu+1} \right).
$$
By Lemma \ref{lemmagejt}, we write
\begin{equation*}
\left\| f_8^{(\lambda,\beta)} \right\|_1 \ge C(\beta)
\sum\limits_{\nu=1}^\infty \lambda_\nu \frac{b_\nu}{\nu}  \ge
C(\beta) \sum\limits_{\nu=0}^\infty \lambda_\nu
\left(\varepsilon_{\nu} - \varepsilon_{\nu+1} \right).
\end{equation*}
This contradicts (\ref{43}). Thus, the series in (\ref{5})
converges.

Let now $p=\infty$. Let $\psi$ be a $Q_{\alpha+\rho,1}$-sequence. Then, we define
\begin{equation*}
f_9(x) = \psi_1 \cos x+ \sum\limits_{\nu=1}^\infty
2^{-\nu(\rho+\alpha)} \left( \psi_{2^\nu} - \psi_{2^{\nu-1}}
\right) \cos  2^\nu x.
\end{equation*}
Similarly, as for the function $f_1$ in the case of $2\le p <\infty,$ it is easy to check that
 $f_9\in H_{\alpha+\rho}^p[\omega]\subset
W_p^{\lambda,\beta}$. Applying Lemma \ref{stlak1}, we get
\begin{eqnarray*}
\left\| f_9^{(\lambda,\beta)} \right\|_\infty &\ge& C(\beta)
\left(\lambda_1 \psi_1+ \sum\limits_{\nu=1}^\infty
\lambda_{2^\nu} 2^{-\nu(\rho+\alpha)} \left(\psi_{2^\nu}
-\psi_{2^{\nu-1}} \right) \right)
\\
& \ge& C(\beta) \sum\limits_{\nu=1}^\infty \lambda_{\nu}
\nu^{-(\rho+\alpha)-1} \psi_{\nu}
\\
& \ge& C(\beta) \sum\limits_{\nu=1}^\infty \left(\lambda_{\nu+1}
- \lambda_{\nu} \right) \omega(1/\nu) =\infty,
\end{eqnarray*}
which contradicts $f_9\in W_p^{\lambda,\beta}$. This completes the proof of the necessity in
(\ref{5}).

{\sf {Step 7.}} We will show the necessity in (\ref{6}) for the case of $\sin \frac{\pi\beta}{2}\ne 0$.
Let $H_{\alpha+r}^p[\omega] \subset W_p^{\lambda,\beta}H_\alpha[\varphi]$.
We take a $q_{\alpha+r, \,1}(\omega)$-sequence $\varepsilon$.
Then (\ref{33}) holds for $\alpha+r$ instead of $\alpha$
and for $1$ instead of $\theta$. Since $\sin \frac{\pi\beta}{2}\ne 0$,
\begin{eqnarray}\label{45}
J&:=& \lambda_{n+1} \omega\left(\frac{1}{n+1}\right) +
n^{-\alpha} \sum\limits_{\nu=1}^n \nu^{r+\alpha}
\left(\nu^{-r}\lambda_{\nu}- (\nu+1)^{-r}\lambda_{\nu+1} \right)
\omega\left(\frac{1}{\nu}\right) \nonumber
\\
&+& | \cos \frac{\beta\pi}{2} | \sum\limits_{\nu=n+2}^\infty
\left(\lambda_{\nu+1}-\lambda_{\nu}\right)
\omega\left(\frac{1}{\nu}\right) + | \sin \frac{\beta\pi}{2} |
\sum\limits_{\nu=n+2}^\infty
 \lambda_\nu
\frac{\omega\left(\frac{1}{\nu}\right)}{\nu} \nonumber
\\
&\le& C(\alpha, \beta, r) \left( \sum\limits_{\nu=n+1}^\infty
 \lambda_\nu \frac{\varepsilon_\nu}{\nu}
+ n^{-\alpha} \sum\limits_{\nu=1}^n \lambda_\nu {\varepsilon_\nu}
\nu^{\alpha-1} \right) =: C(\alpha, \beta, r)
\left(J_1+J_2\right).
\end{eqnarray}
We  will use several times the following evident relations:
\begin{equation}\label{f+-}
\omega_\alpha\Bigl(f, \frac{1}{n}\Bigr)_p \asymp
\omega_\alpha\Bigl(f_{+}, \frac{1}{n}\Bigr)_p +
\omega_\alpha\Bigl(f_{-}, \frac{1}{n}\Bigr)_p,
\qquad\mbox{where}\quad
f_{\pm}(x):= \frac{f(x) \pm f(-x)}{2}.
\end{equation}
\\
{\sf Step 7(р): $p=\infty$ and $\cos \frac{\pi\alpha}{2}\ne 0$}.
In this case by Lemma \ref{my}, we have
$\Bigl($
$f_{7+}^{(\lambda,\beta)}:=\bigl(f_{7}^{(\lambda,\beta)}\bigr)_+$
$\Bigr)$
\begin{eqnarray}\label{vspom78}
\omega_\alpha\Bigl(f_{7}^{(\lambda,\beta)}, \frac{1}{n}\Bigr)_p
&\ge& \omega_\alpha\Bigl(f_{7+}^{(\lambda,\beta)},
\frac{1}{n}\Bigr)_p \ge C (\alpha) n^{-\alpha}
\left\| V_n^{(\alpha)} (f_{7+}^{(\lambda,\beta)})(\cdot) \right\|_p \nonumber \\
 &\ge& C (\alpha) n^{-\alpha} \left|
V_n^{(\alpha)} (f_{7+}^{(\lambda,\beta)})(0) \right|_p
 \ge C (\alpha,\beta) J_2.
 \end{eqnarray}
Using  (\ref{11}) and  $\sum\limits_{k=2n}^{\infty}  a_k    \le  4 E_n(f)_\infty$
 (see \cite{bari-sopr}), we write
\begin{eqnarray}\label{vspom79}
\omega_\alpha\Bigl(f_{7}^{(\lambda,\beta)}, \frac{1}{n}\Bigr)_p
 \ge C (\alpha) E_{[\frac{n}{2}]}(f_{7+}^{(\lambda,\beta)})_p
 \ge C (\alpha,\beta) J_1.
 \end{eqnarray}
 It is proved above that  $f_7\in H_{\alpha+r}^p[\omega] \subset
W_p^{\lambda,\beta}H_\alpha[\varphi].$
Collecting inequalities (\ref{vspom78}), (\ref{vspom79}), and (\ref{45}),
we obtain the estimate in the right-hand side of (\ref{6}).
\\
{\sf Step 7(b): $p=\infty$ and $\cos \frac{\pi\alpha}{2}= 0$}.
If $\cos \frac{\pi\beta}{2}\ne 0,$ then we use (\ref{vspom79}) and
\begin{eqnarray*}
\omega_\alpha\Bigl(f_{7}^{(\lambda,\beta)}, \frac{1}{n}\Bigr)_p
&\ge&
 C (\alpha) n^{-\alpha}
\left\| V_n^{(\alpha)} (f_{7-}^{(\lambda,\beta)})(\cdot)
\right\|_p
 \ge C (\alpha,\beta) J_2.
 \end{eqnarray*}
If $\cos \frac{\pi\beta}{2}=0,$ then $f_7=\pm f_{7+}$ and
\begin{eqnarray*}
\omega_\alpha\Bigl(f_{7}^{(\lambda,\beta)}, \frac{1}{n}\Bigr)_p
 \ge C (\alpha) E_{[\frac{n}{2}]}(f_{7+}^{(\lambda,\beta)})_p
 \ge C (\alpha,\beta) J_1.
 \end{eqnarray*}
 To obtain the estimate of $J_2$, we define
\begin{equation*}
f_{10}(x) = \frac{\varepsilon_0}{2} + \left( \varepsilon_{1} -
\varepsilon_{2}\right) \cos  x + \sum\limits_{\nu=1}^\infty
\left( \varepsilon_{2^{\nu}} - \varepsilon_{2^{\nu+1}} \right)
\cos 2^{\nu} x.
\end{equation*}
It is clear that $E_n(f_{10})_p\le \varepsilon_{n+1}$. Then using (\ref{33}), we get
  $f_{10}\in H_{\alpha+r}^p[\omega] \subset
 W_p^{\lambda,\beta}H_\alpha[\varphi]$.
 Further, applying
$\omega_\alpha\left(f,\frac{1}{n}\right)_p \ge C(\alpha)
n^{-\alpha} \left\|  V_n^{(\alpha)} (f) \right\|_p$
  and Lemma \ref{stlak1}, we write  ($2^m\le n+1 < 2^{m+1}$)
\begin{eqnarray}\label{83}
\omega_\alpha\left(f_{10}^{(\lambda,\beta)},
\frac{1}{2^m}\right)_p &\ge& C(\alpha, \beta) 2^{ - m\alpha}
 \sum\limits_{\nu=0}^{m}
  \lambda_{2^\nu}
2^{\nu\alpha} \left( \varepsilon_{2^{\nu}} -
\varepsilon_{2^{\nu+1}} \right) \nonumber
\\
&\ge& C_1(\alpha, \beta, r) 2^{ - m\alpha}
 \sum\limits_{\nu=0}^{m}
  \lambda_{2^\nu}
2^{\nu\alpha} \varepsilon_{2^{\nu}} - C_2(\alpha, \beta, r)
  \lambda_{2^{m+1}}
\varepsilon_{2^{m+1}}.
\end{eqnarray}
At the same time, $ \omega_\alpha\left(^{(\lambda,\beta)},
\frac{1}{2^m}\right)_p \ge C(\alpha, \beta) \lambda_{2^{m+1}}
\varepsilon_{2^{m+1}}.$
 Then we have
\begin{eqnarray}\label{f10}
\omega_\alpha\left(f_{10}^{(\lambda,\beta)}, \frac{1}{n}\right)_p
\ge C(\alpha, \beta, r) J_2
\end{eqnarray}
 and applying (\ref{f+-}),
\begin{eqnarray*}
C(\alpha, \beta, r) \Bigl[ J_1+J_2\Bigr]  &\le&
 \omega_\alpha\left(f_{10}^{(\lambda,\beta)},
\frac{1}{2^m}\right)_p+
\omega_\alpha\left(f_{7}^{(\lambda,\beta)}, \frac{1}{2^m}\right)_p
\\
&\asymp&
 \omega_\alpha\left(\left(f_{7}+ f_{10}\right)^{(\lambda,\beta)}, \frac{1}{n}\right)_p
= O \left[\varphi\left(\frac{1}{n}\right)\right].
\end{eqnarray*}
The necessity in (\ref{6}) follows.
\\
{\sf Step 7(c): $p=1$ and $\cos \frac{\pi\alpha}{2}\ne 0$}.
In this case we use the function $f_6$.

It is known that $f_6\in H_{\alpha+r}^p[\omega]$ and by Lemma \ref{lemmagejt},
we have
\begin{eqnarray*}
\omega_\alpha\Bigl(f_6^{(\lambda,\beta)}, \frac{1}{n}\Bigr)_1
&\ge& C (\alpha) E_{n}(f_{6-}^{(\lambda,\beta)})_p \ge C
(\alpha,\beta) \sum\limits_{\nu=n+1}^\infty \lambda_{\nu}
\frac{a_\nu}{\nu}
\\
&\ge&C (\alpha,\beta) \left( \sum\limits_{\nu=n+1}^\infty
\lambda_{\nu} \frac{\varepsilon_\nu}{\nu} - a_{n+1}\lambda_{n} -
\sum\limits_{\nu=n+1}^\infty \left(\lambda_{\nu}- \lambda_{\nu-1}
\right) a_\nu  \right).
\end{eqnarray*}
 On the other hand,
\begin{equation*}
\omega_\alpha\Bigl(f_6^{(\lambda,\beta)}, \frac{1}{n}\Bigr)_1 \ge
C (\alpha,\beta,r) \left( a_{n+1}\lambda_{n} +
\sum\limits_{\nu=n+1}^\infty \left(\lambda_{\nu}- \lambda_{\nu-1}
\right) a_\nu \right).
\end{equation*}
Therefore, the last two inequalities imply
$\omega_\alpha\bigl(f_6^{(\lambda,\beta)}, \frac{1}{n}\bigr)_1 \ge
C (\alpha,\beta,r) J_1$. Moreover,
\begin{eqnarray}\label{70}
\omega_\alpha\Bigl(f_{6}^{(\lambda,\beta)}, \frac{1}{n}\Bigr)_p
&\ge&
 C (\alpha) n^{-\alpha}
\left\| V_n^{(\alpha)} (f_{6-}^{(\lambda,\beta)})(\cdot)
\right\|_p
 \ge C (\alpha,\beta) J_2.
 \end{eqnarray}
\\
{\sf Step 7(d): $p=1$ and $\cos \frac{\pi\alpha}{2} =0$}.
 If $\cos \frac{\pi\beta}{2}\ne 0,$ then we use $\omega_\alpha\bigl(f_6^{(\lambda,\beta)},
\frac{1}{n}\bigr)_1 \ge C (\alpha,\beta,r) J_1$ and
\begin{eqnarray*}
\omega_\alpha\Bigl(f_{6}^{(\lambda,\beta)}, \frac{1}{n}\Bigr)_p
&\ge&
 C (\alpha) n^{-\alpha}
\left\| V_n^{(\alpha)} (f_{6+}^{(\lambda,\beta)})(\cdot)
\right\|_p
 \ge C (\alpha,\beta) J_2.
 \end{eqnarray*}
If $\cos \frac{\pi\beta}{2}=0,$ we consider $f_6+f_8$. Using
Lemmas \ref{my} and \ref{lemmagejt}, we get
\begin{eqnarray*}
\omega_\alpha\left(f_8^{(\lambda,\beta)}, \frac{1}{n}\right)_1
&\ge& C(\alpha, \beta) n^{-\alpha}
 \sum\limits_{\nu=1}^{n}
  \lambda_{\nu}
\nu^{\alpha-1} b_\nu
\\
&\ge& C_1(\alpha, \beta, r) n^{-\alpha}
 \sum\limits_{\nu=1}^{n}
  \lambda_{\nu}
\nu^{\alpha-1} \varepsilon_{\nu-1} - C_2(\alpha, \beta, r)
 \lambda_{n}
\varepsilon_{n}.
\end{eqnarray*}
Since $\omega_\alpha\left(f_8^{(\lambda,\beta)},
\frac{1}{n}\right)_1 \ge C(\alpha, \beta, r)
 \lambda_{n}\varepsilon_{n}$, then
$\omega_\alpha\left(f_8^{(\lambda,\beta)},\frac{1}{n}\right)_1
\ge C(\alpha, \beta,r) J_2$.
Thus,
$$
C(\alpha, \beta, r) \Bigl[ J_1+J_2\Bigr] \le
 \omega_\alpha\left(\left(f_{6}+ f_{8}\right)^{(\lambda,\beta)}, \frac{1}{n}\right)_p
=O \left[\varphi\left(\frac{1}{n}\right)\right],$$  i.e.,
 the necessity in (\ref{6}) follows.

{\sf {Step 8.} } We prove the necessity in (\ref{6}) in the case of $\sin \frac{\pi\beta}{2}=
0$.  Let $H_{\alpha+r}^p[\omega] \subset W_p^{\lambda,\beta}H_\alpha[\varphi]$  and let
 $\varepsilon$ be a $q_{\alpha+r, \,1}(\omega)$-sequence.
 Since $\sin \frac{\pi\beta}{2}= 0$, from (\ref{33})  it follows that
\begin{eqnarray}\label{49}
 \lambda_{n+1} \omega\left(\frac{1}{n+1}\right)
&+& n^{-\alpha} \sum\limits_{\nu=1}^n \nu^{r+\alpha}
\left(\nu^{-r}\lambda_{\nu}- (\nu+1)^{-r}\lambda_{\nu+1} \right)
\omega\left(\frac{1}{\nu}\right) \nonumber
\\
&+& | \cos \frac{\beta\pi}{2} | \sum\limits_{\nu=n+2}^\infty
\left(\lambda_{\nu+1}-\lambda_{\nu}\right)
\omega\left(\frac{1}{\nu}\right) \nonumber
\\
&\le& C(\alpha, \beta, r)
 \left(
\sum\limits_{\nu=n+1}^\infty
 \lambda_\nu
\left(\varepsilon_{\nu}-\varepsilon_{\nu+1}\right) + n^{-\alpha}
\sum\limits_{\nu=1}^n \lambda_\nu {\varepsilon_\nu}
\nu^{\alpha-1} \right) \nonumber
\\
&=:& C(\alpha, \beta, r) \left(J_3+J_4\right).
\end{eqnarray}
\\
{\sf Step 8(a): $p=\infty$ and $\cos \frac{\pi\alpha}{2}\ne 0$}.
Applying the Jackson inequality and
Lemma \ref{stlak1}, we get
\begin{eqnarray} \label{f102}
\omega_\alpha\left(f_{10}^{(\lambda,\beta)},
\frac{1}{2^m}\right)_p \ge C(\alpha, \beta)
\sum\limits_{\nu=m}^{\infty}  \lambda_{2^\nu} \left(
\varepsilon_{2^{\nu}} - \varepsilon_{2^{\nu+1}} \right).
\end{eqnarray}
We also note that by Lemma (\ref{stlak1}), (\ref{f10}) holds for all $\alpha>0$. This and
(\ref{f102}) allow us to write
 $\omega_\alpha\left(f_{10}^{(\lambda,\beta)}, \frac{1}{n+1}\right)_p \ge
C(\alpha, \beta) \left(J_3+J_4\right)$.
 Using condition (\ref{49}) and $f_{10}\in H_{\alpha+r}^p[\omega] \subset
W_p^{\lambda,\beta}H_\alpha[\varphi]$,
we obtain the relation in the right-hand side of (\ref{6}).
\\
{\sf Step 8(b): $p=\infty$ and $\cos \frac{\pi\alpha}{2}= 0$}.
Then we consider $f_{10}$ and $f_{11}:=\widetilde{f_{10}}$. It is clear that
$f_{11}\in L_p$ and $f_{10}+f_{11}\in H_{\alpha+r}^p[\omega]$.
 Besides,
\begin{align*}
&\qquad& &\qquad\qquad\qquad
\omega_\alpha\left(f_{10}^{(\lambda,\beta)},
\frac{1}{n+1}\right)_p \ge C(\alpha, \beta) J_3,
 &
\\
&\qquad&
&\qquad\qquad\qquad
\omega_\alpha\left(f_{11}^{(\lambda,\beta)},
\frac{1}{n+1}\right)_p \ge C(\alpha, \beta) J_4
 &
\\
\intertext{
\textnormal{
and}}
&\qquad&
&\omega_\alpha\left(f_{10}^{(\lambda,\beta)},
\frac{1}{n+1}\right)_p+
\omega_\alpha\left(f_{11}^{(\lambda,\beta)},
\frac{1}{n+1}\right)_p \asymp
 \omega_\alpha\left(\left(f_{10}+ f_{11}\right)^{(\lambda,\beta)}, \frac{1}{n+1}\right)_p. &
\end{align*}
{\sf Step 8(c): $p=1$ and $\cos \frac{\pi\alpha}{2}\ne 0$}. Since
$f_8^{(\lambda,\beta)}(x) \sim \pm \sum\limits_{\nu=1}^\infty
 \lambda_\nu b_\nu \sin \nu x,$ then by Lemmas  \ref{my} and \ref{lemmagejt}, we write  (see Step 7(d))
\begin{eqnarray}\label{vspom81}
\omega_\alpha\Bigl(f_{8}^{(\lambda,\beta)}, \frac{1}{n}\Bigr)_p
\ge
 C (\alpha)
 n^{-\alpha}
\left\| V_n^{(\alpha)} (f_{8-}^{(\lambda,\beta)})(\cdot)
\right\|_p
 \ge C (\alpha,\beta) J_4.
 \end{eqnarray}
Further, using  Lemma \ref{lemmagejt} and the Jackson inequality (\ref{11}), we  have
\begin{eqnarray*}
\omega_\alpha\left(f_8^{(\lambda,\beta)},\frac{1}{n}\right)_1
&\ge& C(\alpha, \beta)
 \sum\limits_{\nu=n+1}^{\infty}
  \lambda_{\nu} \frac{b_\nu}{\nu}
\nonumber
\\
&\ge& C(\alpha, \beta)
 \sum\limits_{\nu=n+1}^{\infty}
  \lambda_{\nu}
\sum\limits_{j=2\nu-1}^\infty \frac{\varepsilon_j -
\varepsilon_{j+1}}{j+1} \nonumber
\\
&\ge& C(\alpha, \beta, r) \sum\limits_{j=4n-1}^\infty
\left(\varepsilon_j - \varepsilon_{j+1}\right)
  \lambda_{\nu} .
\end{eqnarray*}
Using the properties of the modulus of smoothness, we get (\ref{6}).
\\
{\sf Step 8(d): $p=1$ and $\cos \frac{\pi\alpha}{2}= 0$}.
We use the fact that  $f_{12}:=\widetilde{f_8} \in L_1$ and
$E_n(f_{12})_p\le C \varepsilon_{n+1}$ (see \cite{gejt}). Then
$f_8+f_{12}\in H_{\alpha+r}^p[\omega]$ and
\begin{eqnarray*}
\omega_\alpha\left(f_8^{(\lambda,\beta)},\frac{1}{n}\right)_1 &\ge&
C(\alpha, \beta) J_3,
\\
\omega_\alpha\left(f_{12}^{(\lambda,\beta)},\frac{1}{n}\right)_1
&\ge& C(\alpha, \beta) J_4.
\end{eqnarray*}

Theorem 1 is fully proved.

\section{Corollaries. Estimates of transformed Fourier series.}
Theorem 1 actually provides estimates of the norms and moduli of smoothness of the transformed
Fourier series, i.e., the estimates of $\|\varphi\|_p$ and $\omega_\alpha(\varphi,\delta)_p$,
where $\varphi \sim \sigma(f, \lambda)$ in terms of $\omega_\gamma(f,\delta)_p$.
 Analyzing the obtained results, one can see that the  following two conditions play a crucial role for these estimates.
 The first is the behavior of the transforming sequence $\{\lambda_n\}$ and the second is the choice between the considered space (as the  Riesz inequality (\ref{9}) holds for $L_p,$ $1<p<\infty$ and no such inequality exists for $L_p,$ $p=1, \infty$).

We will investigate in detail some important examples
 for  $L_p,$ $1<p<\infty$  and  for $L_p,$ $p=1, \infty$, separately.
\\
\begin{center}
{\bf 1. The case of ${\bf 1<p<\infty}$. }
\end{center}
\begin{thr} \label{thr1}
Let\, $1<p<\infty$, $\theta=\min (2,p)$
$\tau=\max (2,p)$, $\alpha\in  \mathbf{R}_+$, $\rho \in
\mathbf{R}_+\cup \{0\}$
 and  $\lambda=\left\{\lambda_n\right\}$ be a
non-decreasing sequence of positive numbers such that  $\left\{n^{-\rho}\lambda_n\right\}$ is non-increasing.
\\
{\bf I. }\,\,If for $f\in L_p^0$ the series
$$
\sum\limits_{n=1}^\infty
\left(\lambda_{n+1}^\theta-\lambda_{n}^\theta\right) \omega^\theta
_{\alpha+\rho}\left(f, \frac{1}{n}\right)_p$$
converges, then there exists a function
$\varphi\in L_p^0$ with the Fourier series $ \sigma(f, \lambda)$, and
\begin{eqnarray}
\label{t1} \|\varphi\|_p &\le& C(p,\lambda, \alpha, \rho) \left\{
\lambda_1^\theta \|f\|_p^\theta + \sum\limits_{n=1}^\infty
\left(\lambda_{n+1}^\theta-\lambda_{n}^\theta\right) \omega^\theta
_{\alpha+\rho}\left(f, \frac{1}{n}\right)_p
\right\}^\frac{1}{\theta},
\\
\omega _{\alpha}\left(\varphi, \frac{1}{n+1}\right)_p &\le&
C(p,\lambda, \alpha, \rho) \left\{ n^{-\alpha\theta}
\sum\limits_{\nu=1}^n \nu^{(\rho+\alpha)\theta}
\left(\nu^{-\rho\theta}\lambda_{\nu}^\theta-
(\nu+1)^{-\rho\theta}\lambda_{\nu+1}^\theta \right) \omega^\theta
_{\alpha+\rho}\left(f, \frac{1}{\nu}\right)_p \right. \nonumber
\\
&+& \sum\limits_{\nu=n+2}^\infty
\left(\lambda_{\nu+1}^\theta-\lambda_{\nu}^\theta\right)
\omega^\theta _{\alpha+\rho}\left(f, \frac{1}{\nu}\right)_p +
\left. \lambda_{n+1}^\theta \omega^\theta _{\alpha+\rho}\left(f,
\frac{1}{n}\right)_p \right\}^\frac{1}{\theta}.\label{t2}
\end{eqnarray}
{\bf II. } \,\,If for $f\in L_p^0$ there exists a function $\varphi\in L_p$ with the Fourier series
$\sigma(f, \lambda)$, then
\begin{eqnarray}
\label{r5} \left\{ \lambda_1^\tau \|f\|_p^\tau +
\sum\limits_{n=1}^\infty \left(\lambda_{n+1}^\tau-\lambda_{n}^\tau
\right) \omega^\tau _{\alpha+\rho}\left(f, \frac{1}{n}\right)_p
\right\}^\frac{1}{\tau} &\le& C(p,\lambda, \alpha, \rho)
\|\varphi\|_p,
\\
\left\{ n^{-\alpha\tau} \sum\limits_{\nu=1}^n
\nu^{(\rho+\alpha)\tau} \left(\nu^{-\rho\tau}\lambda_{\nu}^\tau-
(\nu+1)^{-\rho\tau}\lambda_{\nu+1}^\tau \right) \omega^\tau
_{\alpha+\rho}\left(f, \frac{1}{\nu}\right)_p \right.
&+&\label{r6}
\\
\sum\limits_{\nu=n+2}^\infty \left(
\lambda_{\nu+1}^\tau-\lambda_{\nu}^\tau \right) \omega^\tau
_{\alpha+\rho}\left(f, \frac{1}{\nu}\right)_p +\left.
\lambda_{n+1}^\tau \omega^\tau_{\alpha+\rho}\left(f,
\frac{1}{n}\right)_p \right\}^\frac{1}{\tau} &\le& C(p,\lambda,
\alpha, \rho) \omega_{\alpha}\left(\varphi,
\frac{1}{n+1}\right)_p, \nonumber
\end{eqnarray}
\begin{eqnarray}
\label{tt3} \omega_{\alpha+\rho}\left(f, \frac{1}{n}\right)_p
&\le& C(p,\lambda, \alpha, \rho)
\frac{\|\varphi\|_p}{\lambda_{n}},
\\
\label{tt4} \omega_{\alpha+\rho}\left(f, \frac{1}{n}\right)_p
&\le& C(p,\lambda, \alpha, \rho)
\frac{\omega_{\alpha}\left(\varphi, \frac{1}{n}\right)_p}
{\lambda_{n}}.
\end{eqnarray}
\end{thr}
Inequalities (\ref{t1})-(\ref{t2}) and (\ref{tt3})-(\ref{tt4}) were actually proved in Theorem 1 (see the proof of sufficiency in the part {\bf I}). The estimates (\ref{r5})-(\ref{r6}) are proved analogously, using
Lemma \ref{my}, the theorems by Littlewood-Paley,  Marcinkiewicz, and the Minkowski's inequality (see also \cite{simonov} and \cite{ps1}).

An important corollary of Theorem \ref{thr1} is the following.
\begin{con} \label{con1}
 Let\, $1<p<\infty$, $\theta=\min (2,p)$, and
$\tau=\max (2,p)$. Then for any $k, r>0$ we have
\begin{multline}
\label{b3} C_1 \left\{ \sum\limits_{\nu=n+1}^\infty \nu^{r\tau-1}
\omega_{k+r}^\tau \left(f,\frac{1}{\nu}\right)_p
\right\}^\frac{1}{\tau} \le \omega_{k}
\left(f^{(r)},\frac{1}{n}\right)_p \le  C_2\left\{
\sum\limits_{\nu=n+1}^\infty \nu^{r\theta-1} \omega_{k+r}^\theta
\left(f,\frac{1}{\nu}\right)_p \right\}^\frac{1}{\theta},
\end{multline}
where $C_1=C_1(p,k,r), C_2=C_2(p,k,r), n \in \mathbf{N}.$
\end{con}
See also \cite{studia}, where the right-hand side estimate in (\ref{b3}) was shown for integers  $k$ and $r$.
 The last two inequalities provide sharper bounds in the sense of order than (\ref{b1}) and (\ref{b2}). Indeed, using
 properties of the modulus of smoothness and the Jensen inequality, we have
\begin{eqnarray*}
& &\qquad n^r \omega_{k+r} \left(f,\frac{1}{n}\right)_p \le C(k,r)
 \left\{\sum\limits_{\nu=n+1}^\infty
\nu^{r\tau-1} \omega_{k+r}^\tau \left(f,\frac{1}{\nu}\right)_p
\right\}^\frac{1}{\tau} ;
\\
& & \left\{ \sum\limits_{\nu=n+1}^\infty \nu^{r\theta-1}
\omega_{k+r}^\theta \left(f,\frac{1}{\nu}\right)_p
\right\}^\frac{1}{\theta}
 \le C(k,r) \sum\limits_{\nu=n+1}^\infty
\nu^{r-1} \omega_{k+r} \left(f,\frac{1}{\nu}\right)_p.
\end{eqnarray*}
{\bf Example. } Let $\psi(t)=t^r  \ln^{-A} (1/t)$ and $2\le p
<\infty, \frac{1}{2}<A<1.$ If $\omega_{k+r} \left(f,t\right)_p
\asymp \psi(t)$, then inequalities (\ref{b1}) and (\ref{b2}) give
only $C \ln^{-A} (1/t) \le \omega_{k} \left(f^{(r)}, t
\right)_p$. At the same time, (\ref{b3}) implies $C_1 \ln^{-A+1/p}
(1/t) \le \omega_{k} \left(f^{(r)}, t \right)_p \le C_2
\ln^{-A+1/2} (1/t)$, which is sharper.

Proof of Corollary \ref{con1} follows  from (\ref{t2}) and (\ref{r6}) with $r=\rho$, because if
$f\in L_p,$ $1<p<\infty$, then one can assume that $f^{(r)}\sim \sigma(f, \lambda)$ for $\{\lambda_n=n^r\}$.

The  estimates (\ref{b1}), (\ref{b2}) and  (\ref{b3}) show that it is natural to estimate
$\omega_{\alpha} \left(f^{(\gamma)},\frac{1}{n}\right)_p$
in terms of $\omega_{\alpha+r} \left(f,\frac{1}{\nu}\right)_p$.
 Further analysis allowed us to distinguish three different types of such estimates. It will
be convenient for us to write inequalities in the integral form:

1. $\gamma=r$ (see Corollary \ref{con1})\footnote{Here and further,
$\tau=\max(2,p), \theta=\min(2,p)$. If $A_1\le C A_2$, $C\ge 1$,
we write  $A_1\ll A_2$. Also, if  $A_1\ll A_2$ and $A_2\ll
A_1$, then $A_1\asymp A_2$.}
\begin{equation} \label{weyl1}
\left\{ \int\limits_{0}^{\delta} t^{-r{\tau}-1}
\omega^{\tau}_{r+\alpha}(f,t)_p
 dt \right\}^{\frac{1}{\tau}} \ll \omega_\alpha(f^{(r)},\delta)_p \ll
\left\{ \int\limits_{0}^{\delta} t^{-r{\theta}-1}
\omega^{\theta}_{r+\alpha}(f,t)_p
 dt \right\}^{\frac{1}{\theta}};
\end{equation}

2. $\gamma=r-\varepsilon,$\, $0<\varepsilon<r$ (see Theorem
\ref{thr1} for $\rho=r$ and $\lambda_n=n^{r-\varepsilon}$):
\begin{equation} \label{weyl4}
 \left\{ \int\limits_{0}^{\delta} t^{-(r-\varepsilon){\tau}-1}
 \omega^{\tau}_{r+\alpha}(f,t)_p
 dt + \delta^{\alpha \tau}
\int\limits_{\delta}^{1}  t^{-(r-\varepsilon+\alpha){\tau}-1}
 \omega^{\tau}_{r+\alpha}(f,t)_p dt
 \right\}^{\frac{1}{\tau}}  \ll \omega_\alpha(f^{(r-\varepsilon)},\delta)_p,
\end{equation}
\begin{equation}\label{weyl5}
 \omega_\alpha(f^{(r-\varepsilon)},\delta)_p \ll
\left\{ \int\limits_{0}^{\delta} t^{-(r-\varepsilon){\theta}-1}
 \omega^{\theta}_{r+\alpha}(f,t)_p
 dt + \delta^{\alpha \theta}
\int\limits_{\delta}^{1}  t^{-(r-\varepsilon+\alpha){\theta}-1}
 \omega^{\theta}_{r+\alpha}(f,t)_p dt
 \right\}^{\frac{1}{\theta}};
\end{equation}

3. $\gamma=r+\varepsilon,$\, $0<\varepsilon<\alpha,$ (see
\cite{2003}):
\begin{equation} \label{weyl2}
\left\{ \int\limits_{0}^{\delta} t^{-(r+\varepsilon){\tau}-1}
 \omega^{\tau}_{r+\alpha}(f,t)_p  \,dt\right\}^{\frac{1}{\tau}} \ll
 \delta^{\alpha-\varepsilon}\left\{
\int\limits_{\delta}^{1}  t^{-(\alpha-\varepsilon){\theta}-1}
\omega^{\theta}_\alpha(f^{(r+\varepsilon)},t)_p\, dt
\right\}^{\frac{1}{\theta}},
\end{equation}
\begin{equation} \label{weyl3}
 \delta^{\alpha-\varepsilon}\left\{
\int\limits_{\delta}^{1}  t^{-(\alpha-\varepsilon){\tau}-1}
\omega^{\tau}_\alpha(f^{(r+\varepsilon)},t)_p dt
\right\}^{\frac{1}{\tau}} \ll \left\{ \int\limits_{0}^{\delta}
t^{-(r+\varepsilon){\theta}-1}
 \omega^{\theta}_{r+\alpha}(f,t)_p  dt\right\}^{\frac{1}{\theta}}.
\end{equation}
Some more general estimates of the type (\ref{weyl4})-(\ref{weyl3}) for moduli of smoothness of the transformed Fourier series  can be obtained (\cite{2003}, \cite{2004}, \cite{ps1}) using the sequences of the type
(see, for example, \cite{ulyanov}, \cite{ps1}) $\{\Lambda_n(s):=
\Lambda(s,\frac{1}{n})\}$, where
\begin{equation} \label{posl}
\Lambda(s,t)=\Lambda(s,r,t)=\left( \int\limits^1_t  \xi(u) du +
t^{-rs}\int\limits^t_0 u^{rs} \xi(u) du  \right)^\frac{1}{s}
\end{equation}
and a non-negative function $\xi(u)$ on  $[0,1]$ is such that
$u^{rs} \xi(u)$ is summable.

\begin{center}
{\bf 2. The case of  ${\bf p=1, \infty}$. }
\end{center}

Estimates of $\omega_{\alpha}(\varphi,t)_p$ in terms of $\omega_{r+\alpha}(f,t)_p$
for this case follow from Theorem 1 (see item {\bf II}).
 We will write only the commonly used estimates  of \,$\omega_{\alpha}
(f^{(r)},t)_p$ and $\omega_{\alpha}(\tilde{f}^{(r)},t)_p$ in terms of
$\omega_{r+\alpha}(f,t)_p$ (see also \cite{samko},
\cite{step-1995}-\cite{zhukina}).

\begin{con} \label{con2}
If  $p=1,\infty$, then inequalities (\ref{b1}), (\ref{b2})
hold true for any  $k, r>0.$
\end{con}
If $\left\{\lambda_n=n^\rho\right\}$, $\rho\ge 0$ and
$\beta=\rho+1$, Theorem 1 implies the following
\begin{con} \label{con25}
Let $p=1,\infty$. Then
$$
H_{\alpha+\rho}^p[\omega] \subset \widetilde{W}^\rho_p
\Longleftrightarrow \sum\limits_{n=1}^\infty n^{\rho-1}
\omega\left(\frac{1}{n}\right)\, < \,\infty.
$$
\end{con}
Note that in the case of $p=1, \rho=0$ and $\alpha=1$, Corollary \ref{con25} gives the answer for the question by F.
M\'{o}ricz (\cite{mo}, 1995) on necessary conditions for the embedding $H_{\alpha+\rho}^p[\omega]
\subset \widetilde{W}^\rho_p$. We also mention the papers \cite{bern}, \cite{stechkin}, and \cite{hass},
where the embedding theorems were proved in the necessity part.

\begin{con} \label{con3}
Let $p=1,\infty$ and $r, \alpha, \varepsilon>0$.
\\
{\bf I. }\,\, If for  $f\in L_p$ the series  $\sum\limits_{\nu=1}^\infty
\nu^{r-1} \omega_{r+\alpha+\varepsilon}
\left(f,\frac{1}{\nu}\right)_p$ converges, then there exists
$\tilde{f}^{(r)} \in L_p$ and
\begin{equation*}
\omega_{\alpha} \Bigl(\widetilde{f}^{(r)},\frac{1}{n}\Bigr)_p
\le  C(r,\alpha, \varepsilon) \left( n^{-\alpha}
\sum\limits_{\nu=1}^n \nu^{r+\alpha-1}
\omega_{r+\alpha+\varepsilon}\Bigl(f,\frac{1}{\nu}\Bigr)_p +
\sum\limits_{\nu=n+1}^\infty \nu^{r-1}
\omega_{r+\alpha+\varepsilon} \Bigl(f,\frac{1}{\nu}\Bigr)_p
\right), \,\, n \in \mathbf{N}.
\end{equation*}
{\bf II. } \,\,If for $f\in L_p$ there exists $\tilde{f}^{(r)}\in
L_p$, then
\begin{equation} \label{b5}
\omega_{r+\alpha+\varepsilon} \left(f,\frac{1}{n}\right)_p \le
\frac{C(r, \alpha, \varepsilon)}{n^r} \omega_{\alpha}
\left(\tilde{f}^{(r)},\frac{1}{n}\right)_p,
 \qquad n \in \mathbf{N}.
\end{equation}
\end{con}

Using the direct and inverse approximation theorems, we can
write inequality from the item  {\bf I } of the previous corollary
in the following equivalent form (see also \cite{ba-st}, \cite{samko},
\cite[Vol. 2, Ch. 6 and 7]{step}-\cite{zhukina}):
\begin{con} \label{con35}
Let $p=1,\infty$ and  $r\ge 0, \alpha>0$. If for $f\in L_p$
the series
 $\sum\limits_{\nu=1}^\infty \nu^{r-1} E_\nu\left(f\right)_p$
converges, then there exists $\tilde{f}^{(r)} \in L_p$ and
\begin{equation*}
\omega_{\alpha} \Bigl(\widetilde{f}^{(r)},\frac{1}{n}\Bigr)_p
\le  C(r,\alpha, \varepsilon) \left( n^{-\alpha}
\sum\limits_{\nu=1}^n \nu^{r+\alpha-1}  E_\nu\left(f\right)_p +
\sum\limits_{\nu=n+1}^\infty \nu^{r-1} E_\nu\left(f\right)_p
\right), \qquad n \in \mathbf{N}.
\end{equation*}
\end{con}

\section{Remarks}

1. The embedding theorems for the classes $H^p_l[\omega], {W}^r_p$, and ${W}^r_pH_k[\varphi]$
in the necessity part were investigated, for example,
in the papers by N.K. Bary and S.B. Stechkin (see \cite{ba-st},
$H^p_1[\varphi]\subset H^p_1[\omega]$),
V. {\`E} Ge\u\i t (see \cite{gejt2}, $H^p_l[\varphi]\subset H^p_k[\omega]$,
$H^p_l[\varphi]\subset \widetilde{H}^p_k[\omega]$, $p=1,\infty$),
N. A. Il'yasov (see \cite{il}, $H^p_l[\varphi]\subset W^r_p$,
$H^p_l[\omega]\subset {W}^r_p H_k[\varphi]$, $r,k\in \mathbf{N}$).
Note also that all estimates in Theorem 1 are
 $"$correct$"$ (in the  terminology of Stechkin), that is, the sharpness from the point of view of order, is realized with the help of individual functions.  We thank N. A. Il'yasov for this remark.

Theorem 1 specifies the previous results both in the necessity and
sufficiency parts. Sufficient conditions for the embeddings
$H^p_l[\omega]\subset {W}^{\lambda, \beta}_p H_k[\varphi]$ and
$H^p_l[\omega]\subset {W}^{\lambda, \beta}_p$
were studied in the papers \cite{step}-\cite{zhukina}.
From these articles, particularly, for an important model example $f^{(\lambda,\beta)}\equiv f^{(r)}$, $r>0$,
we have the following estimates
\begin{eqnarray}\label{1thetais12}
E_n(f^{(r)})_p \,&\le&\, C(r) \left(  n^{r}  \, E_n(f)_p + \left\{
\sum\limits_{k=n+1}^\infty k^{r\theta-1} \,
E_k^\theta(f)_p\right\}^\frac{1}{\theta}\right),
\\
\label{1thetais1}
\omega_{k} \Bigl(f^{(r)},\frac{1}{n}\Bigr)_p &\le& \frac{C(k,
r)}{n^k} \left( \sum\limits_{\nu=0}^n (\nu+1)^{(k+r)\theta-1}
E_\nu^\theta(f)_p\right)^\frac{1}{\theta} + C(k, r) \left(
\sum\limits_{\nu=n+1}^\infty \nu^{r\theta-1} E_\nu^\theta(f)_p
\right)^\frac{1}{\theta},
\end{eqnarray}
where $\theta=1$.
At the same time, for the case of $1<p<\infty$, Theorem
\ref{een} and Theorem  \ref{osnth}\footnote{
Theorem 1 implies inequality (\ref{weyl5}), which is equivalent  (using the direct and inverse theorems) to inequality (\ref{1thetais1}).
}
respectively, imply that these estimates hold for $\theta=\min(2,p)$ and this exponent is sharp the best possible.
\\
2.
A generalization of the class ${W}^r_p E[\xi]$ is the class
\begin{equation*}
{W}^{\lambda, \beta}_p E[\omega] = \Bigl\{f\in W_p^{\lambda, \beta}
:\,\,\, E_n \left( f^{(\lambda,\beta)} \right)_p=
O\left[\omega_n\right] \Bigr\}.
\end{equation*}
In this paper we do not consider in detail the embedding theorems between the classes
${W}^{\lambda, \beta}_p E[\omega]$, ${W}^{\lambda,
\beta}_pH_\alpha[\varphi]$, and $E_p[\varepsilon]\equiv {W}^{\{1\},
\,0}_p E[\varepsilon]$.
We only  notice that some results of such types easily follow from direct and inverse theorems
(\ref{a1})-(\ref{a4}), (\ref{asharp}) and some are given in
\cite{og}, \cite[Vol. 2,  Ch. 6 and 7]{step}, \cite{zhukina} and
\cite{h-l}.
For the case when $\lambda_n$ satisfies the $\triangle_2$-condition,
a complete solution of the problem on embedding between the classes $W_p^{
\lambda,\beta}$, ${W}^{\lambda, \beta}_p E[\xi]$, and $E_p[\varepsilon]$
 is described in the following result.
\begin{thr}{\textnormal{\cite{banakh}}}\label{een}\quad
Let \,$1<p<\infty$, $\theta=\min (2,p)$,
$\beta\in \mathbf{R}$, and $\lambda=\left\{\lambda_n\right\}$ be a non-decreasing sequence of positive numbers
 satisfying the  $\triangle_2$-condition, i.e.,  $\lambda_{2n}\le C
\lambda_n$.
Let also $\varepsilon=\{\varepsilon_n\}$ and
$\omega=\{\omega_n\}$ be non-increasing null-sequences.
\\
{\bf I. }\,\, If $1<p<\infty$, then
\begin{eqnarray*}
 E_p[\varepsilon] \subset W_p^{ \lambda,\beta} \quad
&\Longleftrightarrow& \,\,\, \sum\limits_{n=1}^\infty
\left(\lambda_{n+1}^\theta-\lambda_{n}^\theta\right)
\varepsilon^\theta_n \, < \,\infty,
\\
 E_p[\varepsilon] \subset W_p^{ \lambda,\beta}E[\omega] \quad
&\Longleftrightarrow& \biggl\{ \sum\limits_{\nu=n+1}^\infty
\left(\lambda_{\nu+1}^\theta -\lambda_{\nu}^\theta\right)
\varepsilon^\theta_\nu \biggr\}^\frac{1}{\theta} + \lambda_{n}
\varepsilon_n \,\, = O\left[\omega_n \right],
\\
 W_p^{\lambda,\beta} \subset E_p[\varepsilon] \quad
&\Longleftrightarrow& \quad \frac{1}{\lambda_n}= O
\left[\varepsilon_n\right],
\\
 W_p^{ \lambda,\beta}E[\omega]\subset E_p[\varepsilon]
\quad &\Longleftrightarrow& \quad \frac{\omega_n}{\lambda_n}= O
\left[\varepsilon_n\right].
\end{eqnarray*}
{\bf II. } \,\,Let $p=1$ or $p=\infty$.
\\
\textnormal{(a)} If $\triangle\lambda_n \le C
\triangle\lambda_{2n}$ and $\triangle^2\lambda_n\ge 0$
\textnormal{(}or $\le 0$\textnormal{)}, then
\begin{eqnarray*}
E_p[\varepsilon] \subset W_p^{ \lambda,\beta} \quad
&\Longleftrightarrow& \,\,\, | \cos \frac{\beta\pi}{2} |
\sum\limits_{n=1}^\infty  \left(\lambda_{n+1}-\lambda_{n}\right)
\varepsilon_n   \\  &+& \,\,\, | \sin \frac{\beta\pi}{2} |
\sum\limits_{n=1}^\infty \lambda_{n} \frac{\varepsilon_n}{n} \, <
\,\infty,
\\
E_p[\varepsilon] \subset W_p^{ \lambda,\beta}E[\omega] \quad
&\Longleftrightarrow& \,\,\, | \cos \frac{\beta\pi}{2} |
\sum\limits_{\nu=n+1}^\infty
\left(\lambda_{\nu+1}-\lambda_{\nu}\right) \varepsilon_\nu +
\lambda_{n} \varepsilon_n \nonumber
\\
 &+& \,\,\, | \sin \frac{\beta\pi}{2} |
\sum\limits_{\nu=n+1}^\infty \lambda_{\nu}
\frac{\varepsilon_\nu}{\nu} \,\, = O\left[\omega_n \right].
\end{eqnarray*}
\\
\textnormal{(b)}
If for $\beta= 2k, \, k\in \mathbf{Z}$
the condition $\triangle^2\left(1/\lambda_n\right)\ge 0$ holds, and for
$\beta\ne 2k, \, k\in \mathbf{Z}$ the conditions
$\triangle^2\left(1/\lambda_n\right)\ge 0$ and
$\sum\limits_{\nu=n+1}^\infty \frac{1}{\nu\lambda_\nu} \le
\frac{C}{\lambda_n}$  hold, then
\begin{eqnarray*}
W_p^{\lambda,\beta} \subset E_p[\varepsilon] \quad
&\Longleftrightarrow& \quad \frac{1}{\lambda_n}= O
\left[\varepsilon_n\right],
\\
 W_p^{\lambda,\beta}E[\omega]\subset E_p[\varepsilon] \quad
&\Longleftrightarrow& \quad \frac{\omega_n}{\lambda_n}= O
\left[\varepsilon_n\right].
\end{eqnarray*}
\end{thr}
3. The Weyl class $W_p^{\lambda, \beta}$
coincides with the class of functions from $L(0, 2\pi)$
such that their Fourier series can be presented in the following form
$$
\frac{a_0 (f)}{2}+ \sum\limits_{\nu=1}^\infty
 \frac{1}{\pi\lambda_\nu}
\int\limits_0^{2\pi} \psi(x-t) \cos \Bigl(\nu
t-\frac{\pi\beta}{2}\Bigr)\,dt,\quad \psi(x)\in L^0.
$$
Further, consider the case when
$$\sum\limits_{\nu=1}^\infty   \frac{1}{\lambda_\nu}
\cos \Bigl(\nu t-\frac{\pi\beta}{2}\Bigr)$$
is the Fourier series of a summable function $D_{\lambda,\beta}(t)$.

For example, it is certainly so if $\{\lambda_\nu \uparrow \infty\}$ ($n\uparrow$) and
$\sum\limits_{\nu=1}^\infty  \frac{1}{\nu\lambda_\nu}<\infty$
(see \cite[Ch. 5]{zygm}).
Then elements of $W_p^{\lambda, \beta}$ can differ only by the
mean value from functions $f$, which have the following
representation by convolution,
$$f(x)=  \frac{1}{\pi}
\int\limits_0^{2\pi} \psi(x-t)D_{\lambda,\beta}(t) \,dt,\quad
\psi(x)\in L^0.
$$
Here, $\psi$ coincides almost everywhere with $f^{(\lambda,\beta)}$.
See, for example, \cite[Ch. 11]{butzbook}.
Also, similar questions are discussed in detail in the book \cite[Vol. 1, Ch. 3]{step}.
Note that a representation by convolution was first considered in \cite{sz}.
\\
4. The results of sections 2 and 5 can be extended  to a multidimensional case. We only write
 the following estimates for the mixed modulus  of smoothness  $\omega_{\alpha_1,\alpha_2}
(f,\delta_1,\delta_2)_p$ of orders $\alpha_1$ and $\alpha_2$ ($\alpha_1,\alpha_2>0$)
 of a function $f$ (in the $L_p$ metric) with respect to the variables $x_1$ and $x_2$, respectively\footnote{The definition of the mixed modulus  of smoothness and the mixed derivative in the Weyl sense can be found in,
  for example, \cite{2003}.}.
 \begin{thr} {\textnormal{(see also \cite{2003})}} \label{t5}
Let $f(x_1,x_2)\in L_p$, $1<p< \infty$, $\theta= \min(2,p),
\tau=\max(2,p)$ and let
  $\alpha_1, \alpha_2$, $r_1, r_2>0$,
\\
{\bf I.} \,\,If
\begin{equation*}
J_1(\theta):= \left( \int\limits^{1}_{0}\int\limits^{1}_{0}
t_1^{-r_1\theta-1} t_2^{-r_2\theta-1} \omega^\theta_{r_1+\alpha_1,
r_2+\alpha_2}(f,t_1, t_2)_p \,dt_1\, dt_2 \right)^\frac{1}{\theta}
\,<\infty,
\end{equation*}
then $f$ has the mixed derivative in the Weyl sense $f^{(r_1,r_2)}\in L_p^0$.
Moreover,
$$\|f^{(r_1,r_2)}\|_p\le  C(p,r_1,r_2)J_1(\theta)$$  and
\begin{multline*}
\omega_{\alpha_1, \alpha_2} (f^{(r_1,r_2)},\delta_1, \delta_2)_p
\le C(p,\alpha_1, \alpha_2,r_1,r_2) \left(
\int\limits_{0}^{\delta_1}\int\limits_{0}^{\delta_2}
t_1^{-r_1\theta-1} t_2^{-r_2\theta-1} \omega^\theta_{r_1+\alpha_1,
r_2+\alpha_2}(f,t_1, t_2)_p \,dt_1\,  dt_2
\right)^\frac{1}{\theta}
\\=:C\, J_2(\theta).
\end{multline*}
{\bf II.}  \,\,If $f$ has the mixed derivative in the Weyl sense $f^{(r_1,r_2)}\in L_p$,
then
 $$J_1(\tau)\le  C(p,r_1,r_2)  \|f^{(r_1,r_2)}\|_p$$
 and
$$ J_2(\tau)
 \le   C(p,\alpha_1, \alpha_2,r_1,r_2)
 \omega_{\alpha_1, \alpha_2}(f^{(r_1,r_2)},\delta_1, \delta_2)_p.
$$
\end{thr}

\begin{thr} \label{t4}
Let $f(x_1,x_2)\in L_p$, $p=1, \infty$ and let $\alpha_1, \alpha_2, r_1, r_2>0$.
\\
{\bf I.}  \,\, If
\begin{equation*}
J_3:= \int\limits^{1}_{0}\int\limits^{1}_{0} t_1^{-r_1-1}
t_2^{-r_2-1} \omega_{r_1+\alpha_1, r_2+\alpha_2}(f,t_1, t_2)_p
\,dt_1\,  dt_2  \,<\infty,
\end{equation*}
then $f$ has the mixed derivative in the Weyl sense $f^{(r_1,r_2)}\in
L_p$. Moreover,
$$\|f^{(r_1,r_2)}\|_p \le C(r_1,r_2) J_3$$
and
\begin{equation*}
\omega_{\alpha_1, \alpha_2} (f^{(r_1,r_2)},\delta_1, \delta_2)_p
\le C(\alpha_1, \alpha_2,r_1,r_2)
\int\limits_{0}^{\delta_1}\int\limits_{0}^{\delta_2} t_1^{-r_1-1}
t_2^{-r_2-1}
 \omega_{r_1+\alpha_1, r_2+\alpha_2}(f,t_1, t_2)_p \,dt_1\,  dt_2.
\end{equation*}
{\bf II.} \,\,If $f$ has the mixed derivative in the Weyl sense
 $f^{(r_1,r_2)}\in L_p$, then
$$
\omega_{r_1+\alpha_1, r_2+\alpha_2}(f,\delta_1, \delta_2)_p \le
C(\alpha_1, \alpha_2,r_1,r_2) \delta_1^{r_1} \delta_2^{r_2}
\|f^{(r_1,r_2)}\|_p$$ and
\begin{equation*} \label{th02}
\omega_{r_1+\alpha_1, r_2+\alpha_2}(f,\delta_1, \delta_2)_p \le
  C(\alpha_1, \alpha_2,r_1,r_2)
\delta_1^{r_1} \delta_2^{r_2}
 \omega_{\alpha_1, \alpha_2} (f^{(r_1,r_2)},\delta_1, \delta_2)_p.
\end{equation*}
\end{thr}
For more details on the estimates of transformed series in a multidimensional case, see the articles  \cite{2003}-\cite{2004}.
\\
5. In view of inequalities (\ref{weyl1}) and (\ref{weyl4}) -
(\ref{weyl3}), the problem of finding the estimates of
$\omega_{\alpha}\left(\varphi, t\right)_p$ in terms of
$\omega_{\alpha+r}\left(f, t\right)_p$
 arises, e.g., in the case
 $\varphi\sim \sigma(f, \lambda)$, where $\lambda_n=n^r \ln^A
n$. If $A<0$ (which is an analogue of the case $\lambda_n=n^{r-\varepsilon}$),
then estimates $\omega_{\alpha}\left(\varphi, t\right)_p$
follow from \cite{tikh}.

For example, if  $p=2$ and $\varphi\sim \sigma(f, n^r \ln^A n), A<0$,
then
\begin{equation} \label{hhj}
\omega^2_\beta(\varphi,\delta)_2 \asymp
 \int\limits_{0}^{\delta} \frac{t^{-2r-1}}{\ln^{2|A|}\left(\frac{2}{t}\right)}
 \omega^2_{r+\beta}(f,t)_2\, dt + \delta^{2\beta}
\int\limits_{\delta}^{1} \frac{t^{-2(r+\beta)-1}}
{\ln^{1+2|A|}\left(\frac{2}{t}\right)}
 \omega^2_{r+\beta}(f,t)_2 \,dt.
\end{equation}
 Note that the differences between (\ref{hhj}) and (\ref{weyl4})-(\ref{weyl5}) are related
only to the replacement of $n^{-\varepsilon}$ by  $\ln^A n$.

The case $A>0$ (which is an analogue of the case
$\lambda_n=n^{r+\varepsilon}$) is interesting.
For $p=2$ and $\varphi\sim \sigma(f, n^r \ln^A n), A>0$ we have
\begin{equation}\label{hj}
\delta^{2\beta} \ln^{2A}\left(\frac{2}{\delta}\right)
\int\limits_{\delta}^{1} \frac{t^{-2\beta-1}}{
\ln^{1+2A}\left(\frac{2}{t}\right)} \omega^2_\beta(\varphi,t)_2
\, dt +\omega^2_\beta(\varphi,\delta)_2 \asymp
 \int\limits_{0}^{\delta} t^{-2r-1}
\ln^{2A}\left(\frac{2}{t}\right) \omega^2_{r+\beta}(f,t)_2 \, dt.
\end{equation}
Comparing these relations with  estimates (\ref{weyl2})-(\ref{weyl3}), one can remark that the new term
$\omega^2_\beta(\varphi,\delta)_2$ appears in (\ref{hj}).
Thus, this case has essential distinctions. See for detail the papers \cite{2004}, \cite{tikh}.
\\
6. Defining the class $W_p^{\lambda, \beta} H_\alpha[\varphi]$
we assume that $\varphi\in \Phi_\alpha$.
This restriction is natural for  a majorant of the modulus of smoothness of order $\alpha$ (see \cite{tikh-dr}).
\\
7. In Theorem 1 (item II) we used the inequality
$\sum\limits_{\nu=n+1}^\infty \frac{1}{\nu\lambda_\nu} \le
\frac{C}{\lambda_n}$.
It is equivalent to the following condition: there exists  $\varepsilon>0$ such that the sequence
 $\left\{n^{-\varepsilon}\lambda_n\right\}$ is almost increasing, i.e.,
$n^{-\varepsilon}\lambda_n\le C m^{-\varepsilon}\lambda_m$, $C\ge
1$, $n\le m.$
This and other conditions can be found in \cite{ba-st} and \cite{ti}.

The paper was supported by RFFI  (project N 06-01-00268), the programm Leading Scientific Schools
(project NSH-2787.2008.1), Centre de Recerca Matematica and Scuola Normale Superiore.

\vspace{6mm}
\begin{raggedright}
Boris V.  Simonov
\\
Volgograd Technical University
\\ Volgograd, 400131
\\ Russia
\\
e-mail: {\tt htf@vstu.ru}

\vspace{7mm}
Sergey Yu. Tikhonov
\\
ICREA and \\
Centre de Recerca  Matem\`{a}tica
\\ Apartat 50  08193 \, Bellaterra
\\ Barcelona \, Spain
\\
e-mail: {\tt stikhonov@crm.es}\\
\end{raggedright}
\end{document}